\documentclass[10pt]{amsart}
\pdfoutput=1
\usepackage{amsmath}
\usepackage{amssymb}
\usepackage{amsthm}
\usepackage{bm}
\usepackage{accents}
\usepackage{mathtools}
\usepackage{tikz}
\usetikzlibrary{calc}
\usetikzlibrary{decorations.pathmorphing,shapes}
\usetikzlibrary{automata,positioning,quotes}
\usepackage{tikz-cd}
\usepackage{forest}
\usepackage{braket} 
\usepackage{listings}
\usepackage{mdframed}
\usepackage{verbatim}
\usepackage{physics2}
\usephysicsmodule{ab,ab.legacy,diagmat,xmat}
\usepackage{derivative}
\usepackage{fixdif}
\usepackage{stmaryrd}
\usepackage[mathscr]{eucal}
\usepackage{stackengine} 

\usepackage{microtype}

\usepackage{algorithmicx}
\usepackage{algpseudocode}
\usepackage{algorithm}

\usepackage[english]{babel} 
\usepackage[backend=biber,style=alphabetic,maxalphanames=4,maxnames=5,hyperref]{biblatex}
\usepackage[bookmarks, colorlinks, breaklinks]{hyperref} 
\hypersetup{linkcolor=blue,citecolor=magenta,filecolor=black,urlcolor=blue}
\usepackage{cleveref}
\usepackage{xurl}
\usepackage{graphicx}
\graphicspath{{./}}
\crefname{equation}{}{}

\usepackage{float}
\usepackage{booktabs}
\usepackage[shortlabels]{enumitem}
\setitemize{noitemsep}
\usepackage{csquotes}
\newlist{mydescription}{description}{1}
\setlist[mydescription]{style=nextline,
                        font=\bfseries,
                        labelindent=1cm, 
                        leftmargin =2cm,
                        rightmargin=1cm,
                        topsep     =1ex
                       }

\usepackage{lipsum}

\DeclarePairedDelimiter{\floor}{\lfloor}{\rfloor}
\DeclarePairedDelimiter{\ceil}{\lceil}{\rceil}

\newtheorem{thm}{Theorem}[section]
\newtheorem{cor}[thm]{Corollary}

\newtheorem{lem}[thm]{Lemma}

\theoremstyle{definition}
\newtheorem{defn}[thm]{Definition}

\newtheorem{conv}[thm]{Convention}

\theoremstyle{remark}
\newtheorem{rmk}[thm]{Remark}

\theoremstyle{plain}
\newtheorem*{thm*}{Theorem}
\newtheorem*{prop*}{Proposition}
\newtheorem*{lem*}{Lemma}
\newtheorem*{cor*}{Corollary}
\newtheorem*{conj*}{Conjecture}

\theoremstyle{definition}
\newtheorem*{defn*}{Definition}
\newtheorem*{exer*}{Exercise}
\newtheorem*{defns*}{Definitions}
\newtheorem*{con*}{Construction}
\newtheorem*{exm*}{Example}
\newtheorem*{exms*}{Examples}
\newtheorem*{notn*}{Notation}
\newtheorem*{notns*}{Notations}
\newtheorem*{addm*}{Addendum}

\theoremstyle{remark}
\newtheorem*{rmk*}{Remark}

\newcommand{\A}{\mathbb{A}}
\newcommand{\G}{\mathbb{G}}

\newcommand{\R}{\mathbb{R}}
\newcommand{\C}{\mathbb{C}}
\newcommand{\Z}{\mathbb{Z}}
\newcommand{\Q}{\mathbb{Q}}

\newcommand{\ba}{\mathbf{a}}
\newcommand{\bb}{\mathbf{b}}
\newcommand{\bA}{\mathbf{A}}
\newcommand{\bB}{\mathbf{B}}
\newcommand{\bx}{\mathbf{x}}
\newcommand{\bt}{\mathbf{t}}
\newcommand{\bp}{\mathbf{p}}
\newcommand{\bO}{\mathbf{0}}

\renewcommand{\P}{\mathbb{P}}
\newcommand{\E}{\mathbb{E}}

\newcommand{\mc}[1]{\mathcal{#1}}
\newcommand{\msc}[1]{\mathscr{#1}}
\newcommand{\T}{\mathbf{T}}
\newcommand{\mf}[1]{\mathfrak{#1}}
\newcommand{\mbf}[1]{\mathbf{#1}}
\newcommand{\mr}[1]{\mathrm{#1}}
\newcommand{\on}[1]{\operatorname{#1}}
\newcommand{\ms}[1]{\mathsf{#1}}

\newcommand{\ol}[1]{\overline{#1}}

\newcommand{\wh}[1]{\widehat{#1}}

\newcommand{\1}{\mathbf{1}}

\newcommand{\vir}{\mr{vir}}

\newcommand{\pt}{\mr{pt}}
\newcommand{\loc}{\mr{loc}}
\newcommand{\sloop}{\ms{loop}}

\newcommand{\extra}{\ms{extra}}
\newcommand{\tw}{\mr{tw}}
\renewcommand{\top}{\mr{top}}

\DeclareMathOperator{\Aut}{Aut}

\DeclareMathOperator{\ev}{ev}

\DeclareMathOperator{\st}{st}
\DeclareMathOperator{\Res}{Res}

\numberwithin{equation}{section}
\title[Feynman rule and anomaly equations for weighted $\P^4$]{Higher genus Gromov-Witten theory of one-parameter Calabi-Yau threefolds II: Feynman rule and anomaly equations}
\author{Patrick Lei}
\date{\today}

\begin{document}
    
\begin{abstract}
    We prove the Feynman rule conjectured by Bershadsky-Cecotti-Ooguri-Vafa~\cite{bcov} and the anomaly equations conjectured by Yamaguchi-Yau~\cite{yy04} for the Gromov-Witten theory of the Calabi-Yau threefolds $Z_6 \subset \P(1,1,1,1,2)$, $Z_8 \subset \P(1,1,1,1,4)$, and $Z_{10} \subset \P(1,1,1,2,5)$. These determine the generating series $F_g$ of genus $g$ Gromov-Witten invariants recursively from the lower-genus $F_{h<g}$ up to $3g-3$ unknown parameters.
\end{abstract}

\maketitle

\tableofcontents

\section{Introduction}%
\label{sec:Introduction}

\subsection{Historical overview}%
\label{sub:Historical overview}

Despite its origin in theoretical physics as a duality between A-model and B-model topological string theories, mirror symmetry has sparked significant interest in mathematics, starting with the landmark paper~\cite{cdgp}, which gave predictions for the number of genus zero curves of any degree on the quintic threefold. Even though their predictions differ from the actual numbers of curves, they sparked a significant change in the field of enumerative geometry -- namely, to consider deformation-invariant virtual counts of curves. The first mathematical theory constructed to satisfy this property is Gromov-Witten theory (a formalization of A-model invariants), which was developed in symplectic topology by various authors~\cite{rtqcoh,ltsymplectic,siebertsymplecticgw,arnoldgw,ruanvirtualneighborhoods} and in algebraic geometry by Behrend-Fantechi and Li-Tian~\cite{intrinsicnc,gwalggeo,ltgwfoundation}.

For physical reasons, a central problem in Gromov-Witten theory is to compute the Gromov-Witten invariants of compact Calabi-Yau threefolds. For simplicity, we will restrict to the case when $h^2 = 1$ and in particlar to those which arise as complete intersections in weighted projective spaces, of which there are $13$ examples (see~\cite[\S 4]{hkq} for a complete list). We will list (very incompletely) some historical developments in mathematics and physics:
\begin{itemize}
    \item A theorem determining the genus zero invariants generalizing the predictions of~\cite{cdgp} was proved by Givental and Lian-Liu-Yau~\cite{eqgwinv,lly} for complete intersections in projective space and by Coates-Corti-Lee-Tseng and Wang~\cite{orbqcohwproj,mirrornonconvex} for complete intersections in weighted projective spaces.
    \item Bershadsky-Cecotti-Ooguri-Vafa~\cite{bcov} studied the B-model Kodaira-Spencer gravity and predicted that the generating series $F_g$ of genus $g$ Gromov-Witten invariants can be computed recursively using $F_{h<g}$ by a \textit{Feynman rule} up to a finite ambiguity. Mathematically, this corresponds to a sum over stable graphs, which index combinatorial strata of the moduli space $\ol{\msc{M}}_{g,n}$ of stable curves.
    \item Yamaguchi-Yau~\cite{yy04} studied the B-model theory further and predicted that a normalized generating series $P_g$ of genus $g$ Gromov-Witten invariants is a polynomial in five explicit generators for $Z_5 \subset \P^4$, $Z_6 \subset \P(1^4,2)$, $Z_8 \subset \P(1^4,4)$, and $Z_{10} \subset \P(1^3,2,5)$. They also predicted that these polynomials satisfy differential equations in the generators, which we will refer to as \textit{Anomaly Equations}.\footnote{These are typically referred to as Holomorphic Anomaly Equations in the mathematics literature, which is technically inaccurate. In physics, the B-model partition function $F_g^{\bB}(t, \bar{t})$ has an antiholomorphic part, which is governed by the physical Holomorphic Anomaly Equations.} These predictions were extended to the other examples by Huang-Klemm-Quackenbush~\cite{hkq}.
    \item In mathematics, exact formulae for the genus $1$ invariants were proved for complete intersections in projective spaces  by Zinger and Popa~\cite{reducedgenus1,popagenus1ci} and for hypersurfaces in weighted projective spaces by the author~\cite{polynomiality}.
    \item An exact formula for the genus $2$ invariants of the quintic was proved by Guo-Janda-Ruan~\cite{genus2logglsm} pending the proof of foundational results in logarithmic geometry. This was followed by a proof~\cite{bcovlogglsm} of the Anomaly Equations and finite generation conjecture, as well as some other results, for the quintic threefold.
    \item The Feynman rule, Anomaly Equations, and finite generation conjecture were proved independently by Chang-Guo-Li~\cite{nmsp2,nmsp3} using the theory of Mixed-Spin-P (MSP) fields developed in~\cite{mspfermat,msp2,msp3,nmsp}.\footnote{A parameter $N$, which is a positive integer, was introduced in~\cite{nmsp} and the resulting theory was originally referred to as $N$-Mixed-Spin-P (NMSP) fields. In this paper, we will simply refer to the theory with arbitrary $N$ as MSP fields.} Using a slight generalization of their approach, the author proved the finite generation conjecture~\cite{polynomiality} for hypersurfaces in weighted projective space.
\end{itemize}

\subsection{Setup and past work}%
\label{sub:Setup and past work}

Let $\ba = (1,1,1,1,2)$, $(1,1,1,1,4)$, or $(1,1,1,2,5)$, $k \coloneqq \sum_{i=1}^5 a_i$, and set $p_k \coloneqq \frac{k}{a_1\cdots a_5}$. The $I$-function of $Z = Z_k \subset \P(\ba)$ is given by
\begin{align*}
    I(q,z) &\coloneqq z \sum_{d \geq 0} q^d \frac{\prod_{m=1}^{kd} (kH+mz)}{\prod_{i=1}^5 \prod_{m=1}^{a_i d} (a_i H + mz)} \\
    &\eqqcolon I_0(q) z + I_1(q) H + I_2(q) \frac{H^2}{z} + I_3(q) \frac{H^3}{z^2},
\end{align*}
where $H = c_1(\msc{O}_{\P(\ba)}(1))$. If we set $r \coloneqq \frac{k^k}{a_1^{a_1} \cdots a_5^{a_5}}$, then $I(q,z)$ has radius of convergence $\frac{1}{r}$.

\begin{rmk}
    In~\cite{eqgwinv} and other work, the $I$-function differs from ours by a prefactor of $q^{\frac{H}{z}}$. In particular, applying $zq\odv{}{q}$ to the usual conventions corresponds to applying $H+zq\odv{}{z}$ to our $I$-function.

    Our choice is natural from the perspective of quasimap theory. There, our $I$-function is obtained by localization on the stacky loop space~\cite[Proposition 4.9]{orbqmap} and appears in the wall-crossing formula of~\cite{qmapwc}.
\end{rmk}

Define $D \coloneqq q \odv{}{q}$ and define
\[ I_{11} \coloneqq 1 + D\ab(\frac{I_1(q)}{I_0(q)}). \]
Yamaguchi-Yau~\cite{yy04} defined infinitely many generators
\begin{align*}
    A_m \coloneqq \frac{D^m I_{11}}{I_{11}}, \qquad 
    B_m \coloneqq \frac{D^m I_{0}}{I_{0}}, \qquad 
    Y \coloneqq \frac{1}{1-rq}, \qquad \text{and} \qquad X \coloneqq 1-Y.
\end{align*}
For simplicity, denote $A \coloneqq A_1$ and $B \coloneqq B_1$.

\begin{rmk}
    Our generators differ from the generators in~\cite{yy04,hkq} slightly. The variables differ by $\psi^{\ms{HKQ}} = (rq)^{-1}$, and the generators in~\cite{hkq} are defined by
    \[ A_m^{\ms{HKQ}} \coloneqq \frac{\ab(\psi^{\ms{HKQ}} \odv{}{\psi^{\ms{HKQ}}})^m (q I_{11})}{qI_{11}} \qquad \text{and} \qquad B_m^{\ms{HKQ}} \coloneqq \frac{\ab(\psi^{\ms{HKQ}} \odv{}{\psi^{\ms{HKQ}}})^m I_0}{I_0}. \]
    In particular, we have the relations
    \begin{align*}
        A_m = (-1)^m \ab(\sum_{j=0}^m \binom{m}{j} A_j^{\ms{HKQ}}) \qquad \text{and} \qquad
        B_m = (-1)^m B_m^{\ms{HKQ}}.
    \end{align*}
\end{rmk}

\begin{lem}[\cite{yy04}]
    The ring $\msc{R} \coloneqq \Q[A,B,B_2,B_3,Y]$ contains all $A_m$ and $B_m$ for $m \geq 0$. In particular, we have the relations
    \begin{align*}
        A_2 &= 2B^2 - 2AB - 4 B_2 - X (A + 2B + r_0); \\
        B_4 &= -X (2B_3 + (1+r_0)B_2 + r_0 B + r_1),
    \end{align*}
    where $r_0$ and $r_1$ are given in~\Cref{tab:r0andr1}.
    \begin{table}[htpb]
        \centering
        \caption{Values of $r_0$, $r_1$, $a_{0,k}$, and $a_{1,k}$ for different $k$}
        \label{tab:r0andr1}
    
        \begin{tabular}{ccccc}
            \toprule
            $k$ & $r_0$ & $r_1$ & $a_{0,k}$ & $a_{1,k}$ \\ 
            \midrule
            $6$ & $\frac{13}{36}$ & $\frac{5}{162}$ & $\frac{1}{2}$ & $-\frac{7}{4}$ \\
        \\[-0.8em]
            $8$ & $\frac{11}{32}$ & $\frac{105}{4096}$ & $\frac{1}{3}$ & $-\frac{11}{6}$ \\
        \\[-0.8em]
            $10$ & $\frac{3}{10}$ & $\frac{189}{10000}$ & $\frac{1}{6}$ & $-\frac{17}{12}$ \\
            \bottomrule
        \end{tabular}
    \end{table}
\end{lem}

We now define the generating function
\[ F_g(Q) \coloneqq \delta_{g,0} a_{0,k} (\log Q)^3 + \delta_{g,1} a_{1,k} \log Q + \sum_d N_{g,d} Q^d, \]
where the values of $a_{0,k} = \frac{1}{6} \int_Z H^3$ and $a_{1,k} = -\frac{1}{24} \int_Z c_2(Z) \cup H$ are given in~\Cref{tab:r0andr1}.
This is not quite so well-behaved, so we will normalize it by defining
\[ P_{g,m} \coloneqq \frac{(p_k Y)^{g-1} I_{11}^n}{I_0^{2g-2}} \ab(Q \odv{}{Q})^m F_g(Q) \Bigg|_{Q = qe^{\frac{I_1}{I_0}}} \]
for any $(g,m)$ satisfying $2g-2+m > 0$. These satisfy the recursive relation
\[ P_{g,m+1} = (D+(g-1)(2B+X)-mA)P_{g,m} \]
and therefore can be computed from $P_{0,3} = 1$, $P_{1,1}$, and $P_{g \geq 2}$.
The main result of our previous work~\cite{polynomiality} is the following:
\begin{thm}
    For any $(g,m)$ satisfying $2g-2+m > 0$, $P_{g,m} \in \msc{R}$.
\end{thm}

We also proved the following exact formula for the genus $1$ invariants of $Z$:
\begin{thm}
    We have
    \[ P_{1,1} = -\frac{1}{2}A + \ab(\frac{\chi(Z)}{24} - 2) B - \frac{1}{12}X + a_{1,k}. \]
    For clarity, note that $\chi(Z_6) = -204$, $\chi(Z_8) = -296$, and $\chi(Z_{10}) = -288$.
\end{thm}

\subsection{Feynman rule}%
\label{sub:Feynman rule}

\begin{defn}\label{defn:propogators}
    Define the physicists' \textit{propogators} by the formulae
    \begin{align*}
        E_{\psi} &\coloneqq B_1, \\
        E_{\varphi\varphi} &\coloneqq A_1 + 2 B_1, \\
        E_{\varphi\psi} &\coloneqq -B_2, \\
        E_{\psi\psi} &\coloneqq -B_3 + (B_1 - X)B_2 - r_0 B_1 X.
    \end{align*}
    For a choice of ``gauge'' $\G \coloneqq (c_{11}, c_{12}, c_2, c_3)$, where $c_{11},c_{12} \in \Q[X]_{1}$, $c_2 \in \Q[X]_{2}$, and $c_3 \in \Q[X]_{3}$, define
    \begin{align*}
        E_{\psi}^{\G} &\coloneqq E_{\psi} + c_{11}, \\
        E_{\varphi\varphi}^{\G} &\coloneqq E_{\varphi\varphi} + c_{12}, \\
        E_{\varphi\psi}^{\G} &\coloneqq E_{\varphi\psi} - c_{12}B_1 + c_2, \\
        E_{\psi\psi}^{\G} &\coloneqq E_{\psi\psi} + c_{12} B_1^2 - 2 c_2 B_1 + c_3.
    \end{align*}
\end{defn}

Define $I_{22} \coloneqq 1 + D\ab(\frac{D\ab(\frac{I_2}{I_0}) + \frac{I_0}{I_0}}{I_{11}})$ and $I_{33} \coloneqq I_{11}$.
Let $\varphi_i = I_0 \cdots I_{ii} H^i$ for $i = 0,1,2,3$ and let $\psi$ denote the ancestor class on $\ol{\msc{M}}_{g,n}$.
\begin{defn}
    Define the B-model Gromov-Witten correlators $P_{g,m,n}$ by
    \[ P_{g,m,n} \coloneqq \begin{cases}
        (2g+m+n-3)_n P_{g,m} & 2g-2+m > 0 \\
        (n-1)! \ab(\frac{\chi(Z)}{24}-1) & (g,m) = (1,0).
    \end{cases}
    \]
    Here, $(2g+m+n-3)_n$ is the falling Pochhammer symbol.
\end{defn}

Note that these agree with the GW invariants
\[ \frac{(p_k Y)^{m-1}}{I_0^{2g-2+m+n}} \ab<\varphi_1^{\otimes m}, (\varphi_0\psi)^{\otimes n}>_{g,m+n}^Z \]
whenever $(g,m) \neq (1,0)$, in which case the GW invariants is $(n-1)! \frac{\chi(Z)}{24}$. Later, we will discover the meaning of this term.

\begin{defn}\label{defn:bmodelfeynman}
    Let $G_{g,\ell}$ be the set of all stable graphs of genus $g$ and $\ell$ legs and define
    \[ f_{g,m,n}^{\bB,\G} \coloneqq \sum_{\Gamma \in G_{g,m+n}} \frac{\on{Cont}_{\Gamma}^{\bB,\G}}{\ab|\Aut \Gamma|}, \]
    where the contribution of a stable graph is defined by the following construction:
    \begin{itemize}
        \item At each leg, we place $\varphi_1 - E_{\psi}^{\G}\varphi_0 \psi$ or $\varphi_0 \psi$;
        \item At each edge, we place the bivector
            \[ V_{\bB,\G} \coloneqq E_{\varphi\varphi}^{\G} \varphi_1 \otimes\varphi_1 + E_{\varphi\psi}^{\G} (\varphi_1\otimes \varphi_0\psi + \varphi_0\psi\otimes \varphi_1) + E_{\psi\psi}^{\G} \varphi_0\psi\otimes\varphi_0\psi; \]
        \item At each vertex, we place the linear map
            \[ \varphi_1^{\otimes m} \otimes (\varphi_0 \psi)^{\otimes n} \mapsto P_{g,m,n}. \]
    \end{itemize}
\end{defn}

The main result of this paper is the following polynomiality result for $f_{g,m,n}^{\bB}$.
\begin{thm}[\Cref{cor:bmodelfeynman}]\label{thm:bmodelintro}
    For any $g,m,n$ and any choice of gauge $\G$, we have
    \[ f^{\bB,\G}_{g,m,n} \in \Q[X]_{3g-3+m}. \]
\end{thm}

If we specialize to $m=n=0$, then $G_{g,0}$ contains a leading graph $\Gamma_0$ with exactly one vertex (of genus $g$) and no edges. This graph contributes $P_g$, while other graphs contribute linear combinations of products of $P_{h<g,m,n}$ and the propogators. Therefore, if we know all $P_{h<g,m,n}$, the formula
\[ P_g = f_g^{\bB} - \sum_{\Gamma \neq \Gamma_0} \frac{1}{\ab|\Aut \Gamma|} \on{Cont}^{\bB}_{\Gamma} \]
implies that to compute $P_g$, it suffices to compute $f_g^{\bB} \in \Q[X]_{3g-3}$. Because the degree $0$ invariant
\[ N_{g,0} = \frac{(-1)^g \chi(Z) \cdot \ab|B_{2g}| \cdot \ab|B_{2g-2}|}{4g \cdot (2g-2) \cdot (2g-2)!} \]
was already computed by Faber-Pandharipande~\cite{gwhodge}, this allows us to fix the constant term of $f_g^{\bB}$ as $p_k^{g-1} N_{g,0}$.

We will prove this result by reconstructing it from the A-model. We begin by constructing an A-model $R$-matrix, which will act on the Gromov-Witten potential of $Z$. It is engineered such that its edge contributions match the B-model edge contributions as much as possible.
\begin{defn}
    In the basis $\varphi_0, \ldots, \varphi_3$ of $\msc{H}_Z$, define the $A$-model $R$-matrix by
    \[ R^{\bA,\G}(z)^{-1} \coloneqq I - \begin{pmatrix}
        0 & z E_{\psi}^{\G} & z^2 E_{\varphi\psi}^{\G} & z^3 E_{1\psi^2}^{\G} \\
        & 0 & z E_{\varphi\varphi}^{\G} & z^2 E_{1\varphi\psi}^{\G} \\
        & & 0 & z E_{\psi}^{\G} \\
        & & & 0
    \end{pmatrix},
    \]
    where we define $E_{1\varphi\psi}^{\G} \coloneqq -E_{\psi}^{\G} E_{\varphi\varphi}^{\G} - E_{\varphi\psi}^{\G}$ and $E_{1\psi^2}^{\G} \coloneqq -E_{\psi}^{\G} E_{\varphi\psi}^{\G} - E_{\psi\psi}^{\G}$.
\end{defn}

We will now define an A-model Feynman rule.
\begin{defn}
    For $\ba \in \{0,1,2,3\}^n$ and $\bb \in \Z_{\geq 0}^n$, define
    \[ f_{g,m,n}^{\bA,\G} \coloneqq \sum_{\Gamma \in G_{g,m+n}} \frac{\on{Cont}_{\Gamma}^{\bA,\G}}{\ab|\Aut \Gamma|}, \]
    where the contribution of a stable graph is defined by the following construction:
    \begin{itemize}
        \item At each leg, we place $R^{\bA,\G}(z)^{-1} \varphi_a \psi^b$;
        \item At each edge, we place the bivector
            \begin{align*}
                V_{\bA,\G} \coloneqq{}& \sum_{i=0}^3 \frac{\varphi_i \otimes \varphi_{3-i} - R^{\bA}(\psi)^{-1} \varphi_i \otimes R^{\bA}(\psi')^{-1} \varphi_{3-i}}{\psi+\psi'} \\
                ={}& E_{\varphi\varphi} \varphi_1 \otimes \varphi_1 + E_{\varphi\psi} (\varphi_1 \otimes \varphi_0 \psi' + \varphi_0 \psi \otimes \varphi_1) + E_{\psi\psi} (\varphi_0\psi \otimes \varphi_0\psi') \\
                &+ E_{\psi} (\varphi_0 \otimes \varphi_2 + \varphi_2 \otimes \varphi_0) + E_{1\varphi\psi} (\varphi_0 \otimes \varphi_1 \psi' + \varphi_1 \psi \otimes \varphi_0) \\ 
                &+ E_{1\psi^2} (\varphi_0 \otimes \varphi_0 (\psi')^2 + \varphi_0 \psi^2 \otimes \varphi_0).
            \end{align*}
        \item At each vertex, we place the linear map
            \[ \frac{(p_k Y)^{g-1}}{I_0^{2g-2+n}}\ab<->_{g,n}^Z. \]
    \end{itemize}
\end{defn}

Using the theory of MSP fields, we will prove a polynomiality result for $f_{g,\ba,\bb}^{\bA, \G}$. Along the way, we discover that the MSP $R$-matrix defined in~\cite{polynomiality} (see~\Cref{eqn:rmatrix}) factors as the product of a matrix $R^X$ which satisfies an $X$-polynomiality property and $R^{\bA, \G}$ and obtain a similar degree bound for the level $0$ part of the full MSP theory.
\begin{thm}[\Cref{cor:amodelfeynman}]\label{thm:amodelintro}
    For any choice of gauge $\G$ and any $g,n$ such that $2g-2+n > 0$, we have
    \[ f_{g,\ba,\bb}^{\bA, \G} \in \Q[X]_{3g-3+n-\sum b_i}. \]
\end{thm}
In particular, when we set $\ba = (1^m,0^n)$ and $\bb = (0^m,1^n)$, we define
\[ f_{g,m,n}^{\bA, G} \coloneqq f_{g,(1^m,0^n),(0^m,1^n)}^{\bA,\G} \in \Q[X]_{3g-3+m}. \]

Note that the A-model Feynman rule has three extra edge contributions, which contain four extra variables compared to the B-model Feynman rule. Remarkably, the contributions of the correction term in $g=1$, the extra A-model edge contributions, and the four extra A-model variables cancel out to yield the following result.

\begin{thm}[\Cref{thm:allgenusmirror}]\label{thm:aequalsbintro}
    The A-model and B-model graph sums satisfy the relation
    \[ f_{g,m,n}^{\bA, \G} = f_{g,m,n}^{\bB, \G} + \delta_{g,1} \delta_{m,0} (n-1)!. \]
\end{thm}

\subsection{Anomaly equations}%
\label{sub:Anomaly equations}

In~\Cref{sec:Holomorphic anomaly equations}, we prove the following anomaly equations:
\begin{thm}[\Cref{thm:hae}]\label{thm:haeintro}
    The $P_{g,m}$ satisfy the differential equations
    \begin{align*}
    - \partial_{A} P_g = \frac{1}{2} \ab(P_{g-1,2} + \sum_{g_1+g_2 = g} P_{g_1, 1} P_{g_2, 2}),  \\
    \ab(-2 \partial_{A} + \partial_{B} + (A+2B) \partial_{B_2} - \ab((B-X)(A+2B)-B_2-r_0 X)\partial_{B_3})P_g = 0.
    \end{align*}
\end{thm}

The first equation is proved directly by differentiating the B-model Feynman rule, while the second is equivalent to a mysterious reduction of generators (\Cref{thm:reduction}) found by Yamaguchi-Yau~\cite{yy04}. In particular, they consider particular $v_1, v_2, v_3 \in \msc{R}$ and conjecture that for all $g \geq 2$, $P_g \in \Q[v_1, v_2, v_3, X]$.

\subsection{Outline}%
\label{sub:Outline}
The paper is organized as follows:
\begin{itemize}
    \item In \S 2, we review the construction of the MSP $[0,1]$ theory, construct the MSP $[0]$ and $[1]$ theories, and prove polynomiality of the MSP $[1]$ theory.
    \item In \S 3, we prove the polynomiality of the MSP $[0]$ theory and the A-model graph sum using a bootstrapping argument when $\G = (0,0,0,0)$. The A-model Feynman rule for a general gauge (\Cref{thm:amodelintro}) follows as a corollary.
    \item In \S 4, we rewrite both the A-model Feynman rule and B-model Feynman rule using the formalism of geometric quantization of symplectic linear transformations. We then study the effect of the extra contributions in the A-model and prove~\Cref{thm:aequalsbintro} by direct computation. The B-model Feynman rule (\Cref{thm:bmodelintro}) follows as a corollary.
    \item In \S 5, we prove~\Cref{thm:haeintro}.
\end{itemize}

\subsection{Conventions}%
\label{sub:Conventions}

We will use the following conventions in this paper:
\begin{itemize}
    \item We will ignore the odd cohomology of our target $Z$. All operators in this paper will preserve the $\Z/2$ grading on cohomology and are the identity on odd classes.
    \item The theory of MSP fields depends on a positive integer $N$. We will assume that $N$ is an odd prime. Whenever we fix $g,n$, we will always assume that $N \gg 3g-3+n$.
    \item We will consider $\T = (\C^{\times})^N$-equivariant invariants. Our convention is that after equivariant integration, we will specialize our equivariant parameters by $t_{\alpha} = -\zeta_N^{\alpha}t$. At various points in the paper, we will specialize $t^N = -1$.
    \item Whenever we compute using equivariant integration, we will make the substitution $q' = \frac{-q}{t^N}$. Note the specialization $t^N = -1$ makes $q' = q$.
\end{itemize}

\subsection*{Acknowledgements}%
\label{sub:Acknowledgements}

The author is grateful to Chiu-Chu Melissa Liu for all of her helpful advice and for proposing this project. The author would like to thank Konstantin Aleshkin and Shuai Guo for helpful discussions, and Dimitri Zvonkine for his lectures about CohFTs and $R$-matrix actions at the Simons Center for Geometry and Physics in August 2023. The author would also like to thank Shuai Guo for his hospitality during the author's visit to Peking University in July 2024, when part of this work was completed. Finally, the author would like to thank Felix Janda and Yongbin Ruan for expressing interest in the results of this paper and its prequel~\cite{polynomiality}.

\section{MSP $[0]$ and $[1]$ theories}%
\label{sec:MSP 0 and 1 theories}

In this section, we will define the MSP $[0]$ theory and $[1]$ theory and prove a polynomiality property for the $[1]$ theory. We will refer the reader to~\cite{relationsvia3spin} and~\cite[\S 2, Appendix C]{nmsp3} for a discussion of CohFTs and $R$-matrix actions, including in the cases when $R_0 \neq I$ and when source and target of the $R$-matrix are different vector spaces or have different pairings.

\subsection{The MSP $[0,1]$ theory}%
\label{sub:The MSP theory}

Let $N$ be a positive integer. First, the stack $\msc{W}_{g,n,(d,0)}$ was constructed in~\cite{foundations} and MSP invariants were constructed in~\cite{polynomiality} for the state space 
\[ \msc{H} \coloneqq H^*(Z) \oplus \bigoplus_{\alpha=1}^N H^*(\pt_{\alpha}) \eqqcolon \msc{H}_Z \oplus \bigoplus_{\alpha}\msc{H}_{\alpha} \eqqcolon \msc{H}_Z \oplus \msc{H}_1. \]
The pairing is given by 
\begin{align*}
    (x,y) &\coloneqq  \int_{Z} \frac{xy}{-t^N} + \sum_{\alpha} \frac{-p_k}{N t_{\alpha}^3 t^N} xy \bigg\vert_{\pt_{\alpha}}\\
    &\eqqcolon (x|_{Z},y|_{Z})^{Z, \tw} + \sum_{\alpha} (x|_{\pt_{\alpha}}, y|_{\pt_{\alpha}})^{\pt_{\alpha}, \tw}.
\end{align*}
We will now give several bases which we will consider in the rest of the paper.
\begin{enumerate}
    \item We consider $Z \sqcup \bigsqcup_{\alpha=1}^N \pt_{\alpha} = (x_1^{k/a_1} + \cdots + x_5^{k/a_5} = 0)^{\T} \subset \P(\ba, 1^N)$ and let $p = c_1(\mc{O}_{\P(\ba, 1^N)}(1))$. Then define $\phi_j \coloneqq p^j$ for $j = 1, \ldots N+3$;
    \item Note that there is the natural basis $\{1,H,H^2, H^3\}$ of $\msc{H}_Z$ and $\{\1_{\alpha}\}_{\alpha=1}^N$ of $\msc{H}_1 \coloneqq \bigoplus_{\alpha} \msc{H}_{\alpha}$;
    \item We may normalize the previous basis\footnote{This is related to the transformation from flat coordinates to canonical coordinates, see~\cite[\S 5]{polynomiality}} and consider $\varphi_b = I_0 I_{11} \cdots I_{bb} H^b$, where $I_{22} = 1 + D\ab(\frac{D \ab(\frac{I_2}{I_0})+\frac{I_1}{I_0}}{I_{11}})$ and $I_{33} = I_{11}$. We will also consider $\bar{\1}_{\alpha} = L^{-\frac{N+3}{2}}\1_{\alpha}$.
\end{enumerate}

Before we continue, we will define several CohFTs related to $\mc{C}$ using the MSP virtual localization formula~\cite[\S 6]{foundations}. First, for any smooth projective variety $Z$, the Gromov-Witten CohFT associated to $Z$ is given by
\[ \Omega_{g,n}^Z(\tau_1, \ldots, \tau_n) \coloneqq \sum_{d \in H_2(Z,\Z)} q^d \st^Z_* \ab(\prod_{i=1}^n \ev_i^*(\tau_i) \cap [\ol{\msc{M}}_{g,n}(Z,d)]^{\vir}), \]
where $\st^Z \colon \ol{\msc{M}}_{g,n}(Z,d) \to \ol{\msc{M}}_{g,n}$ is the stabilization morphism and $\tau_i \in H^*(Z)$. In the MSP virtual localization formula, the contribution of a vertex at level $0$ is given by
\[ (-t^N)^{-(d+1-g)} [\ol{\msc{M}}_{g,n}(Z,d)]^{\vir} \eqqcolon [\ol{\msc{M}}_{g,n}(Z,d)]^{\tw}. \]
Replacing $[\ol{\msc{M}}_{g,n}(Z,d)]^{\vir}$ by $[\ol{\msc{M}}_{g,n}(Z,d)]^{\tw}$ in the formula for $\Omega^Z$, we obtain the CohFT $\Omega^{Z,\tw}$.

We will need to consider a shift of the Gromov-Witten CohFT of $Z$ by the mirror map $\tau_Z(q) \coloneqq \frac{I_1(q)}{I_0(q)}H$. This is given by the formula
\begin{align*}
    \Omega^{Z,\tau_Z(q)}_{g,n}(-) \coloneqq{}& \sum_{d,m} \frac{q^d}{m!} \Omega^Z_{g,n+m}(-,\tau_Z(q)^m) \\
    ={}& \sum_{d} Q(q)^d \Omega^Z_{g,n+m}(-),
\end{align*}
where $Q(q) \coloneqq qe^{\frac{I_1(q)}{I_0(q)}}$ is the mirror map.

For each of the isolated points $\pt_{\alpha}$, we may consider the vertex contribution
\begin{align*}
    &[\ol{\msc{M}}_{g,n}]^{\alpha, \tw} \\ 
    \coloneqq{} &(-1)^{1-g} \frac{p_k t_{\alpha} \cdot \prod_{i=1}^5 e_{\T}(\E_{g,n}^{\vee}\cdot (-a_i t_{\alpha})) \cdot \prod_{\beta \neq \alpha} e_{\T}(\E_{g,n}^{\vee} \cdot (t_{\beta} - t_{\alpha}))}{(-t_{\alpha})^5 \cdot e_{\T}(\E_{g,n} \cdot kt_{\alpha}) \cdot \prod_{\beta \neq \alpha} (t_{\beta} - t_{\alpha})} \cap [\ol{\msc{M}}_{g,n}].
\end{align*}
These classes define a CohFT $\Omega^{\pt_{\alpha},\tw}$, and restricting to the topological part
\[ [\ol{\msc{M}}_{g,n}]^{\top} = \ab(\frac{1}{p_k} N(-t_{\alpha})^{N+3})^{g-1} [\ol{\msc{M}}_{g,n}] \]
gives the topological part $\omega^{\pt_{\alpha},\tw}$ of the CohFT.

In~\cite[\S 2.1]{polynomiality}, we defined MSP invariants $\ab<->^M_{g,n}$ using virtual localization, whose explicit formula was proved in~\cite[\S 5]{foundations}. Recall that for any cohomological field theory, Givental's theory~\cite{symplfrob} considers the fundamental solution of the quantum differential equation (or Dubrovin connection), which is given by
\[ S^M_{\tau}(z)x = x + \sum_{a, n} \frac{1}{n!} e^a \ab<\frac{x}{z-\psi}, e_a, \tau^n>_{0,n+2}^M. \]
Here, $\{e_a\}$ is a basis of $\msc{H}$ and $\{e^a\}$ is its dual basis. We will now restrict to the case $\tau = 0$ and abbreviate the fundamental solution to $S^M(z)$. We also considered the corresponding fundamental solutions $S^Z \coloneqq S^Z_{\tau_Z(q)}$ where $\tau_Z = \frac{I_1(q)}{I_0(q)}H$ and $S^{\pt_{\alpha}} \coloneqq S_{\tau_{\alpha}}^{\pt_{\alpha}}$, where $\tau_{\alpha} = -t_{\alpha} \int_0^{q} (L(x)-1) \frac{\d{x}}{x}$, where $L = (1+rx)^{\frac{1}{N}}$.

Finally, we defined the MSP $R$-matrix by the formula
\begin{equation}\label{eqn:rmatrix}
    S^M(z) \begin{pmatrix}
    \on{diag}\ab\{ \Delta^{\pt_{\alpha}}(z)\}_{\alpha=1}^N & \\
    & 1
\end{pmatrix} = R(z) \begin{pmatrix}
    \on{diag}\ab\{ S^{\pt_{\alpha}}(z)\}_{\alpha=1}^N & \\
    & S^{Z}(z)
\end{pmatrix}\biggr|_{q \mapsto q'},
\end{equation}
where $\Delta^{\pt_{\alpha}}(z)$ is the Quantum Riemann-Roch~\cite{qrr} operator given by the formula
\begin{align*}
    \Delta^{\pt_{\alpha}}(z) \coloneqq \exp &\Biggl[\sum_{m \geq 0} \frac{B_{2m}}{2m(2m-1)} \Biggl(\sum_{i=1}^5 \frac{1}{(-a_i t_{\alpha})^{2m-1}}  \\
    &+ \frac{1}{(kt_{\alpha})^{2m-1}} + \sum_{\beta \neq \alpha} \frac{1}{(t_{\beta} - t_{\alpha})^{2m-1}}\Biggr) z^{2m-1} \Biggr].
\end{align*}
Here, the $B_{2k}$ are the Bernoulli numbers.

In~\cite[\S 5]{polynomiality}, we explained how to use the explicit formula for the quantum differential equation given in~\cite[Lemma 2.18]{polynomiality} to compute the entries of $R(z)^*$ when the input basis is $\{1, p, \ldots, p^{N+3}\}$ and the output basis is $\{\varphi_0, \ldots, \varphi_3\} \cup \{\bar{\1}_{\alpha}\}_{\alpha=1}^N$. In particular, the entries are elements of $\msc{R}$ up to some normalization.

\begin{defn}[{\cite[Theorem 3.6]{polynomiality}}]
    Define the local theory by the formula
    \[ \Omega^{\loc} \coloneqq \Omega^{Z,\tw} \oplus \bigoplus_{\alpha=1}^N \omega^{\pt_{\alpha},\tw} \]
    and the MSP $[0,1]$ theory by the formula
    \[ \Omega^{[0,1]} \coloneqq R.\Omega^{\loc}. \]
\end{defn}

\begin{rmk}
    Note that the $[0,1]$ theory was originally defined using virtual localization. The fixed loci of $\mc{W}_{g,n,(d,0)}$ are described using localization graphs $\Theta$ whose vertices can be partitioned as $V = V_0 \sqcup V_1 \sqcup V_{\infty}$. We then only consider those $\Theta$ for which $V_{\infty} = \emptyset$ when defining the $[0,1]$ theory. In~\cite[\S 3]{polynomiality}, we proved that it is equivalent to the definition we give here. In addition, when $N$ is very large relative to $g,n$, we obtain a simpler formula for the tail contributions at level $0$.
\end{rmk}

\begin{rmk}
    In contrast to the usual setting (see~\cite{relationsvia3spin} for example), our $R(z) = R_0 + R_1 z + \cdots$ does not satisfy $R_0 = \mr{Id}$. However, we can relate this more general case of $R$-matrix actions to the usual setting via the dilaton flow, for example see~\cite[Appendix C]{nmsp3}.
\end{rmk}

\subsection{The MSP $[0]$ and $[1]$ theories}%
\label{sub:The MSP 0 and 1 theories}

We will now define restricted versions of the MSP $[0,1]$ CohFTs, which we will call the $[0]$ and $[1]$ theories. First, recall that~\cite[\S 2]{nmsp3} gives a definition of the $R$-matrix action on CohFTs when the source and target of $R$ are not the same vector space. In particular, we only need that $R(-z)^* R(z) = \mr{Id}$, which in particular implies that $R_0$ is injective.

\begin{defn}
    Define the restricted $R$-matrices $R^{[0]}(z)$ and $R^{[1]}(z)$ by the formulae
    \begin{align*}
        R^{[0]}(z) &= R(z)|_{\msc{H}_Z}, \\
        R^{[1]}(z) &= R(z)|_{\msc{H}_1}.
    \end{align*}
\end{defn}

Because the MSP $R$-matrix $R(z)$ satisfies $R(-z)^* R(z) = \mr{Id}$, $R^{[0]}(z)$ and $R^{[1]}(z)$ satisfy the identities
\begin{align*}
    R^{[0]}(-z)^* R^{[0]}(z) &= \mr{Id}_{\msc{H}_Z}, \\
    R^{[1]}(-z)^* R^{[1]}(z) &= \mr{Id}_{\msc{H}_1}, \\
    R^{[0]}(-z)^* R^{[1]}(z) &= R^{[1]}(-z)^* R^{[0]}(z) = 0 .
\end{align*}

\begin{defn}
    Define the MSP $[0]$ theory by the formula
    \[ \Omega^{[0]} \coloneqq R^{[0]}. \Omega^{Z,\tw} \]
    and the MSP $[1]$ theory by the formula
    \[ \Omega^{[1]} \coloneqq R^{[1]}. \bigoplus_{\alpha=1}^N \omega^{\pt_{\alpha},\tw}. \]
\end{defn}

Our immediate goal is now to prove a polynomiality result for the MSP $[0]$ theory. We will do this by studying the MSP $[0,1]$ theory and the $[1]$ theory, which is similar to the argument used to prove the polynomiality of the $[0,1]$ theory~\cite[Theorem 4.1]{polynomiality}. We will first describe a bipartite graph decomposition of the $[0,1]$ theory, then prove the polynomiality of the $[1]$ theory. After some work, we will apply the polynomiality of the $[1]$ theory and of the $[0,1]$ theory to deduce the polynomiality of the $[0]$ theory.

\begin{defn}
    Define $\msc{G}^{[0,1]}_{g,n}$ to be the set of stable bipartite graphs of total genus $g$ and $n$ legs. These are stable graphs with a partition $V = V_0 \sqcup V_1$ making the graph bipartite.
\end{defn}

\begin{thm}\label{thm:01bipartite}
    There is a decomposition of the MSP $[0,1]$ theory in terms of stable bipartite graphs as
    \begin{align*}
        \Omega^{[0,1]}_{g,n}(\tau_1, \ldots, \tau_n) = \sum_{\Lambda \in \msc{G}_{g,n}^{[0,1]}} & \bigotimes_{v \in V_0} \Omega^{[0]}_{g_v, n_v} \otimes \bigotimes_{v \in V_1} \Omega^{[1]}_{g_v, n_v} \\
        &\ab(\bigotimes_{\substack{v \in V_0 \\ \ell \in L_v}} \tau_{\ell} \otimes \bigotimes_{\substack{v \in V_0 \\ \ell \in L_v}} \tau_{\ell} \otimes \bigotimes_{e \in E} V^{01}(\psi,\psi')),
    \end{align*}
    where we define
    \[V^{01}(z,w) \coloneqq \sum_{\alpha=1}^N \frac{R^{[1]}(z) - R^{[1]}(-w)}{z+w} \1_{\alpha} \otimes R^{[1]}(w) \1^{\alpha}. \]
\end{thm}

\begin{proof}
    Recall that the edge contribution in the definition of the MSP $[0,1]$ theory is given by
    \begin{align*}
        \on{Cont}_{E_{01}} &= \sum_{i=1}^{N+3} \frac{\phi_i \otimes \phi^i - R(z)^{-1} \phi_i \otimes R(w)^{-1}\phi^i}{z+w} \Bigg\vert_{\msc{H}_Z \otimes \msc{H}_1} \\
        &= - \sum_{\alpha=1}^N \frac{R(-z)^* R(-w) \1_{\alpha} \otimes \1^{\alpha}}{z+w} \Bigg\vert_{\msc{H}_Z \otimes \msc{H}_1} \\
        &= - \sum_{\alpha=1}^N \frac{R^{[0]}(-z)^* R^{[1]}(-w) \1_{\alpha} \otimes \1^{\alpha}}{z+w} \\ 
        &= \sum_{\alpha=1}^N \frac{R^{[0]}(-z)^* (R^{[1]}(z) - R^{[1]}(-w))\1_{\alpha \otimes \1^{\alpha}}}{z+w} \\
        &= (R^{[0]}(-z)^* \otimes R^{[1]}(-w)^* V^{01}(z,w)).
    \end{align*}
    The result then follows from the definition of the $R$-matrix action as a sum over stable graphs.
\end{proof}

\subsection{Polynomiality of the $[1]$ theory}%
\label{sub:Polynomiality of the 1 theory}

Recall that we defined the edge contributions to the $[0,1]$ theory by
\begin{align*}
    V(z,w) &= \sum_{j=0}^{N+3} \frac{\phi_j \otimes \phi^j - R(z)^{-1}\phi_j \otimes R(w)^{-1}\phi^j}{z+w} \\ 
    &\eqqcolon \sum_{m,n} V_{mn} z^m w^n.
\end{align*}

\begin{lem}\label{lem:rlevel1}
    Let $m,n \geq 0$, $a = 0,\ldots, N+3$, and $\alpha,\beta \in \{1,\ldots,N\}$. In~\cite[\S 5.3]{polynomiality}, we defined
    \[ (R_m)_a^{\alpha} \coloneqq L_{\alpha}^{-a-m} (R_m \bar{\1}^{\alpha},\phi_a). \] Then
    \begin{enumerate}
        \item $(R_m)_a^{\alpha}$ is independent of $\alpha$ and $(R_m)_a^{\alpha} \in \Q[X]_{m+\floor{\frac{a}{N}}}$;
        \item The $V$-bivector has the form
            \[ V_{mn} |_{\msc{H}_1 \otimes \msc{H}_1} = L^{-3} t^N \sum_{\alpha,\beta} \sum_s L_{\alpha}^{s-m} L_{\beta}^{2-s-n} (V_{mn})^{\alpha\beta;s} \1_{\alpha} \otimes \1_{\beta}, \]
            where $(V_{mn})^{\alpha\beta;s} \in \Q[X]_{m+n+1}$ is independent of $\alpha,\beta$.
    \end{enumerate}
\end{lem}

\begin{proof}
    Note that (1) is~\cite[Lemma 5.9]{polynomiality}. The proof of (2) is the same as~\cite[Lemma C.1]{nmsp2}.
\end{proof}

\begin{defn}
    When $\star$ is either $[0]$, $[1]$, or $[0,1]$, define
    \[ f^{\star}_{g,(\ba,\bb)} \coloneqq \int_{\ol{\msc{M}}_{g,n}} \prod_{i=1}^n \psi_i^{b_i} \Omega_{g,n}^{\star}(\phi_{a_1}, \ldots, \phi_{a_n}) \]
    and the $[1]$ theory with special insertions
    \begin{align*}
        f^{[1]}_{g,(\ba,\bb)(\ba',\bb')} \coloneqq L^{\sum_{i=1}^m a_i'}& \int_{\ol{\msc{M}}_{g,n+m}} \prod_{i=1}^n \psi_i^{b_i} \prod_{j=1}^m \psi_{n+j}^{b_j'} \\
        & \Omega_{g,n+m}^{[1]}(\phi_{a_1}, \ldots, \phi_{a_n}, R(\psi_{n+1}) \bar{\phi}_{a_1'},\ldots,R(\psi_{n+m}) \bar{\phi}_{a_m'})
    \end{align*}
    for $\ba \in \{0,\ldots,N+3\}^n$, $\ba' \in \{1,\ldots,N\}^m$, $\bb \in \Z_{\geq 0}^n$, and $\bb' \in \Z_{\geq 0}^m$. Here, we define $\bar{\phi}_a \coloneqq L^{-\frac{N+3}{2}} L^a p^a|_{\msc{H}_1} = \sum_{\alpha=1}^N L_{\alpha}^a \bar{\1}_{\alpha}$.
\end{defn}

\begin{defn}
    For a tuple $\ba \in \Z_{\geq 0}^n$, define $\ab|\ba| \coloneqq \sum_{i=1}^n a_i$ and $\floor*{\frac{\ba}{N}} \coloneqq \sum_{i=1}^n \floor*{\frac{a_i}{N}}$.
\end{defn}

\begin{lem}\label{lem:polynomiality1theory}
    Let $N \gg 3g-3+n+m$. Define
    \[ c \coloneqq \frac{\ab|\ba| + \ab|\ba'| + \ab|\bb| + \ab|\bb'| - n-m}{N}. \]
    If $c \in \Z$, then
    \[ \ab( \frac{Y}{t^N} )^{g-1+r} f_{g,(\ba,\bb),(\ba',\bb')}^{[1]} \]
    is a polynomial in $X$ of degree at most $3g-3+n+m-\ab|\bb|-\ab|\bb'| + \floor*{\frac{\ba}{N}}$. Otherwise, $f_{g,(\ba,\bb),(\ba',\bb')}^{[1]} = 0$.
\end{lem}

\begin{proof}
    Recall that $f^{[1]}_{g,(\ba,\bb),(\ba',\bb')}$ is defined as a sum of stable graph contributions, where the contribution from any stable graph is computed as follows:
    \begin{enumerate}
        \item At every leg with insertion $\phi_a \psi^b$, we place
            \[ R^{[1]}(-\psi)^* \psi_a \psi^b = \sum_{\alpha,m} L_{\alpha}^{a-m}(-1)^m \psi^{m+b} \bar{\1}_{\alpha} =  L^{-\frac{N+3}{2}}\sum_{\alpha,m} L_{\alpha}^{a-m}(-1)^m \psi^{m+b} \1_{\alpha}; \]
        \item At every leg with special insertion $R(\psi')\bar{\phi}_{a'}\psi^{b'}$, we place
            \[ R^{[1]}(-\psi)^* R(\psi) \bar{\phi}_{a'} \psi^{b'} = \sum_{\alpha} L_{\alpha}^{a'} \bar{\1}_{\alpha} \psi^{b'} = L^{-\frac{N+3}{2}} \sum_{\alpha} L_{\alpha}^{a'} \1_{\alpha}\psi^{b'}; \]
        \item At every edge, we place the bivector\footnote{Note that we reindex the sum from~\Cref{lem:rlevel1}.}
            \[ V(z,w)|_{\msc{H}_1 \otimes \msc{H}_1} = L^{-3} t^N \sum_{\alpha,\beta} \sum_{c,d,s} L_{\alpha}^{1+s-c} L_{\beta}^{1-s-d} (V_{cd})^{\alpha\beta;s+1} \1_{\alpha} \otimes \1_{\beta}; \]
        \item At every vertex of genus $g_v$ with $n_v$ legs, we place the map
            \[ \sum_{\alpha} L^{\frac{N+3}{2}(2g_v - 2 + n_v)} \sum_s \frac{1}{s!} \st^s_* \omega_{g_v,n_v+s}^{\pt_{\alpha},\tw} (-,T_{\alpha}(\psi)^s), \]
            where
            \begin{align*}
                T_{\alpha}(z) &= z(\1- L^{\frac{N+3}{2}} R(z)^{-1} \1) |_{\pt_{\alpha}} \\
                &= \sum_{m=1}^{\infty} L_{\alpha}^{-m} (R_m)_0^{\alpha} (-z)^{m+1} \1_{\alpha}
            \end{align*}
            was defined in~\cite[Theorem 3.6]{polynomiality} and $\st^s \colon \ol{\msc{M}}_{g,n+s} \to \ol{\msc{M}}_{g,n}$ is the morphism forgetting the last $s$ marked points. Recall from~\cite[Lemma 2.4]{polynomiality} that
            \[ \omega_{g,n}^{\pt_{\alpha},\tw} = \ab(\frac{1}{p_k}N(-t_{\alpha})^{N+3})^{g-1} \Omega_{g,n}^{\pt}. \]
    \end{enumerate}

    We will now count the degrees of the contributions at each vertex labeled by $\pt_{\alpha}$. Let $L_v$ be the set of ordinary (first $n$) legs attached to $v$ and $L_v'$ be the set of special legs attached to $v$. The factor involving $L_{\alpha}$, $L$, and $Y$ is
    \begin{align*}
        L_{\alpha}^{(N+3)(g_v-1)} \prod_{\ell \in L_v} L_{\alpha}^{a_{\ell}-c_{\ell}} \prod_{\ell' \in L_v'} L_{\alpha}^{a'_{\ell'}} \prod_{\substack{f =(e,v)\\e \in E_v}} t^{\frac{N}{2}} L^{\frac{N}{2}} L_{\alpha}^{1+s_f-c_f} \prod_{\ell'' \in \mr{tails}} L_{\alpha}^{-c_{\ell''}}
    \end{align*}
    Using the fact that
    \begin{align*}
        \sum_{\ell' \in \mr{tails}} c_{\ell''} + \sum_{\substack{f = (e,v) \\ e \in E_v}} c_f + \sum_{\ell \in L_v} c_{\ell} + b_{\ell} + \sum_{\ell' \in L_v'} b'_{\ell'} &= 3g_v - 3 + n_v  \\
        &= 3g_v -3 + \ab|L_v| + \ab|L_v'| + \ab|E_v|,
    \end{align*}
    the contribution becomes
    \begin{align*}
        & (tL)^{N\ab(g-1 + \frac{\ab|E_v|}{2})} L_{\alpha}^{3g-3+n_v + \sum_{\ell} (a_{\ell}-1) + \sum_{\ell'} (a'_{\ell'}-1) + \sum_f s_f - \sum_{\ell} c_{\ell} -\sum_f c_f - \sum_{\ell''} c_{\ell''}} \\
        ={}& (tL)^{N\ab(g-1+\frac{\ab|E_v|}{2})} L_{\alpha}^{\sum_{\ell} (a_{\ell}+b_{\ell}-1) + \sum_{\ell'} (a'_{\ell'} + b'_{\ell'}-1) + \sum_f s_f}.
    \end{align*}
    
    The remaining tail, edge, and leg contributions contribute a degree in $X$ of at most
    \begin{align*}
        & \sum_{\ell \in L_v} \ab(c_{\ell} + \floor*{\frac{a_{\ell}}{N}}) + \sum_{\substack{f = (e,v) \\ e \in E_v}} \ab(c_f + \frac{1}{2}) + \sum_{\ell'' \in \mr{tails}} c_{\ell''} \\
        ={}& 3g_v - 3 + n_v + \frac{\ab|E_v|}{2} + \sum_{\ell \in L_v} \floor*{\frac{a_{\ell}}{N}} - \sum_{\ell \in L_v} b_{\ell} - \sum_{\ell' \in L_v'} b'_{\ell'}
    \end{align*}
    
    Note that the end result must be independent of the hours $\alpha$, so we can sum over all $\alpha$ and obtain a multiplicative factor
    \[ L_{\alpha}^{\frac{1}{N}\ab(\sum_{\ell \in L_v}(a_{\ell}+b_{\ell}-1) + \sum_{\ell' \in L_v'} (a'_{\ell'} + b'_{\ell'}-1) + \sum_{e \in E_v} s_{(e,v)})}. \]
    Denote the exponent by $c_v$. Clearly if this is not an integer, the contribution vanishes after summing over all $\alpha$. By our choice of indexing of the edge contributions, we see that $s_{(e,v_1)} + s_{(e,v_2)} = 0$ whenever $e = (v_1, v_2)$. In particular, taking the product over all vertices gives an exponent of
    \[ \sum_v c_v = \frac{1}{N} (\ab|\ba| + \ab|\ba'| + \ab|\bb| + \ab|\bb'| - n-m) = c. \]
    If this is not an integer, the total contribution vanishes.

    Multiplying the prefactors over all vertices, we obtain a total prefactor
    \[ L_{\alpha}^{N \sum_v r_v} (tL)^{\sum_v N \ab(g_v-1+\frac{\ab|E_v|}{2})} = (tL)^{Nc + N(g-1)}. \]
    Using the fact that $Y = L^{-N}$, the prefactor becomes
    \[ \ab(\frac{Y}{t^N})^{-(g-1+c)}. \]

    Multiplying the remaining contributions from the tails, edges, and legs, the total degree in $X$ is at most
    \begin{align*}
        & \sum_v 3g_v - 3 + n_v + \frac{\ab|E_v|}{2} + \sum_{\ell \in L_v} \floor*{\frac{a_{\ell}}{N}} - \sum_{\ell \in L_v} b_{\ell} - \sum_{\ell' \in L_v'} b'_{\ell'} \\
        ={}& 3g-3+n+m + \floor*{\frac{\ba}{N}} - \ab|\bb| - \ab|\bb'|.
    \end{align*}
    The desired result follows after applying the normalization factor $\ab(\frac{Y}{t^N})^{g-1+c}$.
\end{proof}

\subsection{Vanishing of the $[0]$ theory}%
\label{sub:Vanishing of the 0 theory}

We will now study the MSP $[0]$ theory. Unfortunately, we need to factorize the MSP $[0]$ theory to prove the desired polynomiality result, but we will prove a vanishing result similar to the first part of~\Cref{lem:polynomiality1theory}.

\begin{defn}
    Define the \textit{mod-$N$ degrees} by $\deg \varphi_j = \deg \phi_j = j \mod{N}$ and $\deg \psi = 1$.
\end{defn}

\begin{lem}\label{lem:RpreservesmodNdegree}
    $R^{[0]}(z)$ preserves the mod-$N$ degree. If we define $\bar{j} = j\mod{N}$, then
    \[ R^{[0]}(z)^* \phi_j = c'_{j,k} q^{\floor*{\frac{j}{N}}} \varphi_{\bar{j}} + O(z^{\bar{j}-3}) \]
    for $j = 0,\ldots,N+3$, where we define\footnote{These are in fact the $q$ coefficients of the quantities $I_0$, $I_0 I_{11}$, $I_0 I_{11} I_{22}$, and $I_0^2 I_{11} I_{22}$.}
    \begin{align*}
        c'_{j,6} &= (1,\ldots,1,-360,-3132,-8532,-11304) \\
        c'_{j,8} &= (1,\ldots,1,-1680,-17488,-48048,-63856) \\
        c'_{j,10} &= (1,\ldots,1,-15120,-194640,-605360,-784880).
    \end{align*}
\end{lem}

\begin{proof}
    Recall that
    \[ zD R^{[0]}(z)^* = R^{[0]}(z)^* A^M - A^Z R^{[0]}(z)^*, \]
    where
    \[ A^Z = \begin{pmatrix}
        0 \\
        I_{11} & 0 \\
        & I_{22} & 0 \\
        & & I_{11} & 0
    \end{pmatrix}
    \]
    and $A^M$ is given by the formulae
    \[ A^M_{j+1,j} = 1, \qquad A^M_{j+N-1,j} = c_{j,k}q - \delta_{i,4}t^N \]
    and all other entries being zero, where
    \begin{align*}
        c'_{j,6} &= (360,2772,5400,2772,360) \\
        c'_{j,8} &= (1680,15808,30560,15808,1680) \\
        c'_{j,10} &= (15120,179520,410720,179520,15120). 
    \end{align*}
    The expression for $R^{[0]}\phi_j$ follows directly. The fact that $R^{[0]}$ preserves the mod-$N$ degree follows from the fact that
    \[ R^{[0]}(z)^{-1}x = S^Z(q',z) (S^M(z)^{-1}x)|_Z \]
    for all $x \in \msc{H}_Z$ and the fact that both $S^Z$ and $S^M$ preserve the mod-$N$ degrees. Here, it is helpful to recall that
    \begin{align*}
      S^{Z}(z)^* =&{}\ I + \frac{1}{z} \begin{pmatrix}
        0 \\
        J_1' & 0 \\
               & J_2' & 0 \\
               & & J_1' & 0
      \end{pmatrix} \\ &+ \frac{1}{z^2} \begin{pmatrix}
        0 \\
        & 0 \\
        J_2 & & 0 \\
          & \frac{J_2'}{J_1'}J_1 - J_2 & &  0
      \end{pmatrix} + \frac{1}{z^3} \begin{pmatrix}
        0 \\
        & 0 \\
        & & 0 \\
        J_3 & & & 0
      \end{pmatrix},
    \end{align*}
    where we define $J_b = \frac{I_b}{I_0}$, $J_1' = I_{11}$, and $J_2' = J_1 + D J_2$.
\end{proof}

\begin{lem}
    Define $c \coloneqq \frac{1}{N}(\ab|\ba| + \ab|\bb| - n)$. If $c \notin \Z$, then $f_{g,(\ba,\bb)}^{[0]} = 0$.
\end{lem}

\begin{proof}
    Recall that $\Omega^{[0]} = R^{[0]}.\Omega^{Z,\tw}$ is defined by a graph sum formula, where the vertex, edge and leg contributions are given by the following:
    \begin{itemize}
        \item At each leg with insertion $\phi_a \psi^b$, we place $R^{[0]}(-\psi)^* \phi_a \psi^b$;
        \item At each edge, we place the bivector
            \[ \sum_{j=0}^{N+3} \frac{\phi_j \otimes \phi^j - R^{[0]}(-z)^* \phi_j \otimes R^{[0]}(-w)^* \phi^j}{z+w}; \]
        \item At each vertex, we place the linear map $I_0(q')^{-(2g-2+n)} \Omega_{g,n}^{Z,\tw,\tau_Z(q')}(-)$.
    \end{itemize}
    Because $\phi_j \otimes \phi^j$ has mod-$N$ degree $3$ and $R^{[0]}$ preserves the mod-$N$ degree, we see that the edge contributions have mod-$N$ degree $2$. Then note that for $\bullet$ being either ``$Z$'' or ``$Z,\tw$'' the integral
    \[ \int_{\ol{\msc{M}}_{g,m}} \Omega_{g,m}^{\bullet}\ab(\bigotimes_{j=1}^m \varphi_{a_j} \psi^{b_j}) \]
    vanishes unless $\sum_{j=1}^m (a_j+b_j) = m$. Summing over all vertices and edges in an arbitrary stable graph, we see that the contribution can only be nonzero if
    \[ \sum_{i=1}^n (\bar{a}_i + \bar{b}_i) = n, \]
    which is equivalent to $c \in \Z$.
\end{proof}

We may now assume that $c \coloneqq \frac{1}{N}(\ab|\ba| + \ab|\bb| - n)$ is an integer.

\begin{lem}\label{lem:vanishing0nospecial}
    Define
    \[ \bar{c} \coloneqq \frac{\ab|\bar{\ba}| + \ab|\bb| - n}{N}. \]
    If $N \gg 3g-3+3n$, we always have $\bar{c} \geq 0$. In addition, if $\bar{c} \neq 0$, then $f_{g,(\ba,\bb)}^{[0]} = 0$.
\end{lem}

\begin{proof}
    Using the assumption that $N > 3g-3+3n$ and the fact that $\ol{\msc{M}}_{g,n}$ is a DM stack (hence $3g-3+n \geq 0$), we obtain $N > 2n$. Because all $a_i,b_i \geq 0$ and $\bar{c}$ is an integer, we must have
    \[ \bar{c} \geq \ceil*{\frac{-n}{N}} = 0. \]

    The vanishing result is another degree-counting argument. Assume that $\bar{c} > 0$. Recall $f_{g,(\ba,\bb)}^{[0]}$ is a sum of stable graph contributions, where we place
    \[ R^{[0]}(-\psi)^* \phi_{a_i} \psi^{b_i} \]
    at the $i$-th leg. In particular, the total degree of the ancestors is at least
    \begin{align*}
        \ab|\bar{\ba}| - 3n + \ab|\bb| &= N\bar{c} - (\ab|\bb| - n) - 3n + b \\
        &\geq N - 2n.
    \end{align*}
    However, the contribution vanishes if
    \[ \ab|\bar{\ba}_v| - 3(n_v - \ab|E_v|) + \ab|\bb_v| > 3g_v - 3 + n_v \]
    and therefore if
    \begin{align*}
        \ab|\bar{\ba}| - 3n + \ab|\bb| &> \sum_v (3g_v - 3 + n_v) \\
        &= 3g - 3 + n - \ab|E|. 
    \end{align*}
    However, by assumption, we see that
    \begin{align*}
        \ab|\bar{\ba}| - 3n + \ab|\bb| &\geq N - 2n \\
        &> 3g-3+n \\
        &\geq 3g-3+n-\ab|E|,
    \end{align*}
    so all graph contributions must vanish.
\end{proof}

\begin{cor}\label{cor:criterionforc}
    If $f_{g,(\ba,\bb)}^{[0]}$ is nonzero, then $c = \floor*{\frac{\ba}{N}}$ and
    \[ g-1+c \leq 3g-3+c+n-\ab|\bb|. \]
\end{cor}

\begin{defn}\label{defn:0potspecial}
    Define the MSP $[0]$ potential with special insertions by
    \[ f_{g,(\ba,\bb),(\ba',\bb')}^{[0]} \coloneqq \int_{\ol{\msc{M}}_{g,n}} \Omega_{g,n+m}^{[0]}\ab(\bigotimes_{i=1}^n \phi_{a_i} \psi^{b_i}, \bigotimes_{j=1}^m E_{a'_j, b'_j}(\psi_{n+j})), \]
    where we define
    \[ E_{a',b'}(\psi) = L^{-a'} \cdot [w^{b'}] \frac{( R(\psi) - R(-w) )\bar{\phi}^{a'}}{\psi+w}. \]
    Here, $[w^{b'}]$ means that we take the coefficient of $w^{b'}$ in the expression.
\end{defn}

By construction, we have
\[ V^{01}(z,w) = \sum_{a=1}^N E_{ab}(z) w^b \otimes L^a R(w) \bar{\phi}_a. \]

\begin{lem}
    If $a-m \not\equiv b \pmod{N}$, then $(\phi_a, R_m \bar{\phi}^b) = 0$. In the case when $a-m \equiv b \pmod{N}$, then
    \[ L^{-a+m}(\phi_a, R_m \bar{\phi}^b) \in \Q(t^N)[X]_{m+\floor*{\frac{a}{N}}}. \]
\end{lem}

\begin{proof}
    Note that
    \begin{align*}
        (\phi_a, R_m \bar{\phi}^b) &= \sum_{\alpha=1}^N L_{\alpha}^{-b} (\phi_a, R_m \bar{\1}^{\alpha}) \\
        &= \sum_{\alpha=1}^N L_{\alpha}^{-b} L^{a-m} (R_m)_a^{\alpha}.
    \end{align*}
    The result follows by using~\Cref{lem:rlevel1} and the the fact that the total power of the roots of unity is $a-m-b$.
\end{proof}

\begin{lem}
    The $[0]$ potential with special insertions $f_{g,(\ba,\bb),(\ba',\bb')}^{[0]}$ vanishes unless
    $c \coloneqq \frac{\ab|\ba| + \ab|\bb| + \ab|\ba'| + \ab|\bb'| -n+m}{N} \in \Z$. 
\end{lem}

\begin{proof}
    Write
    \[ R_m \bar{\phi}^{b} = \sum_{s=1}^{N+3} (R_m \bar{\phi}^b, \phi_s) \phi^s. \]
    By the previous lemma, the only nonzero terms are the ones with $s \equiv m+b (\pmod N)$, so the mod-$N$ degree of $R_m\bar{\phi}^b$ is $3-(m+b)$. This implies that the mod-$N$ degree of $E_{ab}$ is $2-(a+b)$. Computing the total ancestor degree at all vertices as in the proof of~\Cref{lem:vanishing0nospecial}, we obtain the desired result.
\end{proof}

\section{The $A$-model Feynman rule}%
\label{sec:The A model Feynman rule}

In order to extract information from the $[0]$ theory, we will extract a part of it which satisfies an $X$-polynomiality property.

\subsection{Factorization of the $[0]$ theory}%
\label{sub:Factorization of the 0 theory}

\begin{defn}
    Define the CohFT\footnote{Note that the unit of this CohFT is $\varphi_0$ and not $1$.} $\Omega^{\bA, \G}$ by the formula
    \[ \Omega^{\bA, \G} = R^{\bA,\G}.\Omega^{Z,\tw}. \]
    Also, define the generating function
    \[ f_{g,(\ba,\bb)}^{\bA,\G} \coloneqq \ab( \frac{-p_kY}{t^N} )^{g-1} \int_{\ol{\msc{M}}_{g,n}} \Omega_{g,n}^{\bA,\G}(\varphi_{a_1} \psi_1^{b_1}, \ldots, \varphi_{a_n} \psi_n^{b_n}). \]
\end{defn}

\begin{rmk}
    Because the CohFT $\Omega^Z$ and $\Omega^{\bA}$ have different units, the graph sum formula for the $R^{\bA}$ action will have $T = z(1-I_0)$, which by the dilaton equation will imply that we place the linear map
    \[ \frac{1}{I_0^{2g-2+n}} \Omega^Z_{g,n} \]
    at each vertex.
\end{rmk}

\begin{rmk}
    Because the basis $\varphi_0, \ldots, \varphi_3$ is not flat, we will need to consider the transformation
    \[ \Psi \coloneqq \pdiagmat[empty={}]{I_0,I_0 I_{11},I_0 I_{11} I_{22},I_0I_{11^2}I_{22}} \]
    from this basis to the flat basis $1, H, H^2, H^3$. In particular, if we want to compute derivatives of $R^{\bA,\G}$, then we will need the input basis to be the flat basis, and in this basis $R^{\bA, \G}(z)^{-1}$ takes the form
    \[ \Psi \ab(I - \begin{pmatrix}
        0 & z E_{\psi}^{\G} & z^2 E_{\varphi\psi}^{\G} & z^3 E_{1\psi^2}^{\G} \\
        & 0 & z E_{\varphi\varphi}^{\G} & z^2 E_{1\varphi\psi}^{\G} \\
        & & 0 & z E_{\psi}^{\G} \\
        & & & 0
    \end{pmatrix}). \]
\end{rmk}

\begin{defn}
    Define the matrix
    \[ R^{X}(z) \coloneqq R^{\bA}(z)^{-1} \cdot R^{[0]}(z). \]
    We will consider the input basis of $R^X$ to be $\varphi_0, \ldots, \varphi_3$ and the output basis to be $\phi_0, \ldots, \phi_{N+3}$.
\end{defn}

\begin{lem}\label{lem:matrixelementsRX}
    The matrix elements $(R^X_m)_j^a \coloneqq (\phi_j, R_m^X \varphi^a)$ of $R^{X}$ satisfy the following whenever $m < N-3$:
    \begin{enumerate}
        \item If $j \not\equiv m+a \pmod{N}$, then $(R^X_m)_j^a = 0$;
        \item If $j < N$, then $(R^X_m)_j^a \in X \Q[X]_{m-1}$;
        \item If $j \geq N$, then $(R^X_m)_j^a \in q \Q[X]_m$;
        \item We have $R^X(-z)^* C^X(z) = \mr{Id}_{\msc{H}_Z}$ for some $C^X(z)$ of the form $\begin{psmallmatrix}
            C(z) \\ 0
        \end{psmallmatrix}$, where $C(z)$ has the form
            \[ C(z) =  \begin{pmatrix}
                1 & z \cdot C_{1} & z^2 \cdot C_{2} & z^3 \cdot C_{3}\\
                0 & 1 & z \cdot C_{4} & z^2 \cdot C_{5} \\
                0 & 0 & 1 & z \cdot C_6 \\
                0 & 0 & 0 & 1
            \end{pmatrix}
            \]
            for some $C_1, C_2, C_4, C_6 \in \Q[X]_1$ and $C_3, C_5 \in \Q[X]_2$.
    \end{enumerate}
\end{lem}

\begin{proof}
    Using the equation
    \begin{equation}\label{eqn:mspqde}
     zD R^{[0]}(z)^* = R^{[0]}(z)^* A^M - A^Z R^{[0]}(z)^* . 
    \end{equation}
    Using the definition $R^{[0]}(z)^* = R^{\bA}(z)^* R^X(z)^*$, we obtain
    \[ zD(R^{\bA}(z)^*R^X(z)^*) = R^{\bA}(z)^* R^X(z)^* A^M - A^Z R^{\bA}(z)^* R^X(z)^*, \]
    and multiplying by $R^{\bA}(-z)$, we obtain
    \begin{equation}\label{eqn:Xqde}
    R^X(z)^* A^M = zD(R^X(z)^*) + R^{\bA}(-z) [(zD+A^Z) R^{\bA}(z)^*] R^X(z)^*. 
    \end{equation}
    Direct computation (here, note that the basis $\varphi_0, \ldots, \varphi_3$ is not flat, so we need the version of $R^{\bA}$ with the $\Psi$) yields
    \[ R^{\bA}(-z) [(zD+A^Z)R^{\bA}(z)^*] = \begin{pmatrix}
        0 & 0 & 0 & -r_1 X z^4 \\
        1 & 0 & -r_0 X z^2 & 0 \\
        0 & 1 & -Xz & 0 \\
        0 & 0 & 1 & -Xz
    \end{pmatrix}.
     \]
     Therefore, we can compute $R^X(z)^* \phi_i$ from $R^X(z)^* \phi_0$. Because $R^{[0]}(z)^* \phi_0 = \varphi_0 + O(z^{N-3})$ and $R^{\bA}(z)^* \varphi_0 = \varphi_0$, the first entry of $R^X$ is $1 + O(z^{N-3})$.

     We then note that the matrices $A^M$ and $R^{\bA}(-z)[(zD+A^Z)R^{\bA}(z)^*]$ both increase the mod-$N$ degree by $1$. Therefore, we see that $R^{X}$ preserves the mod-$N$ degree, so the vanishing property holds. The degree estimates follow from the explicit formulae for $A^M$ and $R^{\bA}(-z) [(zD+A^Z)R^{\bA}(z)^*]$. The final statement is obtained by direct computation.
\end{proof}

\begin{cor}\label{cor:edgeX}
    The edge contribution
    \[ V_X \coloneqq \frac{\sum_{i=0}^3 \varphi_i \otimes \varphi^i - \sum_{j=0}^{N+3} R^X(-z)^* \phi_j \otimes R^X(-w)^* \phi^j}{z+w} \]
    satisfies the degree bound
    \[ \frac{Y}{t^N} [z^{m_1} w^{m_2}] V_X \in \msc{H}_Z^{\otimes 2}[X]_{m_1+m_2+1}. \]
\end{cor}

\begin{proof}
    This follows from the lemma and the fact that $\varphi_j = \frac{t^N}{p_k Y} \varphi_{3-j}$.
\end{proof}

\subsection{Polynomiality of the $[0]$ theory and the $\bA$ theory}%
\label{sub:Polynomiality of the 0 theory}

We now prove the polynomiality of the $[0]$-theory. First, we will prove some lemmas about $E_{a',b'}(z)$, which was introduced in~\Cref{defn:0potspecial}.

\begin{lem}\label{lem:vanishingcoeffedge}
    Recall the definition of $E_{a',b'}(z)$ from~\Cref{defn:0potspecial}. Then 
    \[ (\phi^a, [z^b] E_{a',b'}(z)) = 0 \]
    unless $a+a'+b+b' \equiv 2 \pmod{N}$.
\end{lem}

\begin{proof}
    By definition, we have
    \begin{align*}
        (\phi_a, [z^b] E_{a',b'}(z)) &= (-1)^{b'}(\phi_a, R_{b+b'+1} \bar{\phi}^{a'}) \\
        &= (-1)^{b'}  \sum_{\alpha} L_{\alpha}^{-a'} (\phi_a, R_{b+b'+1} \bar{\1}^{\alpha}),
    \end{align*}
    which vanishes unless $a-a'-(b+b'+1) \equiv 0 \pmod{N}$. The result follows from the fact that $\phi^a$ has mod-$N$ degree $3-a$.
\end{proof}

\begin{lem}\label{lem:degreecoeffedge}
    Whenever $N \gg b''$, then
    \[ \ab(\frac{Y}{t^N})^{-c_E} \cdot [z^{b''}] (\varphi^{a''}, R^X(-z)^* E_{a',b'}(z)) \in \Q[X]_{b'+b''+1}, \]
    where we define $c_E \coloneqq \frac{a'+b'+a''+b''-N-2}{N}$. If $c_E \notin \Z$, then the above quantity vanishes.
\end{lem}

\begin{proof}
    The vanishing is a corollary of~\Cref{lem:matrixelementsRX} and~\Cref{lem:vanishingcoeffedge}. The degree estimate comes from the following considerations.

    Whenever $a = 4, \ldots, N-1$ and $m' < N-3$, then $(\varphi^{a''}, [z^{m'}] R^X(-z)^* \phi_a)$ is nonzero only if $a = a''+m'$. If this is satisfied, then it has degree $m'$ in $X$, and therefore
    \[ \ab(\frac{Y}{t^N})^{-c_E} \cdot [z^{b''}] (\varphi^{a''}, R^X(-z)^* \phi_a)(\phi^a, E_{a',b'}(z)) \in \Q[X]_{b'+b''+1}. \]
    Here, we use the fact that 
    \[ [z^b](\phi^a, E_{a',b'}(z)) = \frac{N}{p_k} \ab(\frac{Y}{t^N})^{c_E} (-1)^{b'} (R_{b+b'+1})_{N+3-a}^{\alpha} \in \ab(\frac{Y}{t^N})^{c_E} \Q[X]_{b+b'+1} \] by~\Cref{lem:rlevel1} because $\floor*{\frac{a}{N}} = 0$.

    In the case where $a = 0,\ldots,3$, then for any $m' < N-3$, the vanishing still holds, so we use $\phi^a = \frac{1}{5}(\phi_{N+3-a} - t^N \phi_{3-a})$ to compute. The $\phi_{N+3-a}$ term contributes an element of $\Q[X]_{b'+b''+2}$, whereas the $t^N \phi_{3-a}$ term contributes
    \[ [z^b](t^N \phi_{3-a}, E_{a',b'}(z)) = N (-1)^{b'} \ab(\frac{Y}{t^N})^{c_E} (R_{b+b'+1})_{3-a}^{\alpha} Y. \]
    Therefore, we have
    \[ \ab(\frac{Y}{t^N})^{-c_E} \cdot [z^{b''}] (\varphi^{a''}, R^X(-z)^* \phi_a)(\phi^a, E_{a',b'}(z)) \in \Q[X]_{b'+b''+2}. \]

    In the case when $a \geq N$, the nonvanishing condition becomes $a-N = a''+m'$. Therefore,~\Cref{lem:matrixelementsRX} implies that
    \[ \ab(\frac{Y}{t^N})^{-c_E} \cdot [z^{b''}] (\varphi^{a''}, R^X(-z)^* \phi_a)(\phi^a, E_{a',b'}(z)) \in q\Q[X]_{b'+b''+1}. \]

    If we sum the contributions from the above procedure, the degree count is too high by $1$. The $X^{b'+b''+2}$-coefficient is given by (up to a constant)
    \begin{align*}
        [X^{b'+b''+2}] \sum_{\substack{i+j=b'+b''+1 \\ j \leq b'' \\ a \leq 3}} (R^X_j)_i^a (R_i)_{N+3-a} - Y (R_j^X)_i^a(R_i)_{3-a} + \frac{Y}{t^N} (R_j^X)_{i}^{N+a} (R_i)_{3-a}.
    \end{align*}
    Note that the MSP quantum connection~\Cref{eqn:mspqde} (in particular the value of $A^M$) implies that $[X^{m+1}] (R_m)_j = \frac{c'_{j,k}}{r} \cdot [X^m] (R_m)_{j-N}$ and the equation~\Cref{eqn:Xqde} for $R^X$ implies that $[X^m] (q^{-1} (R^X_m)_j^a) = c'_{j,k} \cdot [X^m] (R_m)_{j-N}^a$, where $c'_{j,k}$ was defined in~\Cref{lem:RpreservesmodNdegree}. This implies that
    \begin{align*}
        & [X^{i+j+1}] ( (R_j^X)_i^a (R_i)_{N+3-a} - Y(R_j^X)_i^a (R_i)_{3-a} + \frac{Y}{t^N} (R_j^X)_i^{N+a}(R_i)_{3-a}) \\
        ={}& \ab(1+\frac{c'_{N+a}}{r} + \frac{c'_{N+3-a}}{r}) \\
        ={}& 0
    \end{align*}
    for all $i+j = b'+b''+1$ and $a = 0,1,2,3$, where we have used the fact that $c'_{N+a} + c'_{N+3-a} = -r$, which is implied by the fact that $I_0^2 I_{11}^2 I_{22} = Y$.
\end{proof}

We will now put a partial ordering on the set of pairs $(g,n)$ such that $(h,m) \prec (g,n)$ if $(h,m) < (g,n)$ in the lexicographic order and $3h+m \leq 3g+n$. We will use induction on $(h,m)$ under this ordering and a bootstrapping argument to prove both the polynomiality of the $[0]$ theory and of the $\bA$ theory.

We introduce the following statements:
\begin{enumerate}
    \item Denote by $\mf{P}_{g,n}$ the statement ``for all $\ba \in \{0,1,2,3\}^n$ and $\bb \in \Z_{\geq 0}^n$, we have 
        \[ \ab(\frac{Y}{t^N})^{g-1} f_{g,(\ba,\bb)}^{[0]} \in \Q[X]_{3g-3+n-\ab|\bb|}.\text{''} \]
    \item Denote by $\mf{Q}_{g,s}$ the statement ``for all $m+n \leq s$, $\ba \in \{0,\ldots,N+3\}^m$, $(\bb,\bb') \in \Z_{\geq 0}^{m+n}$, and $\ba' \in \{1,\ldots,N\}^n$, 
        \[\ab(\frac{Y}{t^N})^{g-1} f_{g,(\ba,\bb),(\ba,\bb')}^{[0]} \in \Q[X]_{3g-3+m+2n+\floor*{\frac{\ba}{N}}-\ab|\bb|+\ab|\bb'|}.\text{''} \]
\end{enumerate}

\begin{lem}\label{lem:Apolyassuming0poly}
    Suppose that $\mf{P}_{h,m}$ holds for all $(h,m) \prec (g,n)$. Then for all $(h,m) \prec (g,n)$, we have
    \[ f_{h,(\ba,\bb)}^{\bA} \in \Q[X]_{3h-3+m-\ab|\bb|} \]
    for all $\ba \in \{0,1,2,3\}^m$.
\end{lem}

\begin{proof}
    First, note that $R^{\A}$ preserves degrees, so we must have $\sum_i (a_i + b_i) = n$. Now define
    \[ \tilde{f}_{h,(\ba,\bb)}^{[0]} \coloneqq \int_{\ol{\msc{M}}_{h,m}} (R^X.\Omega^{\bA})_{h,m}(C^X(\psi_1)\varphi_{a_1}\psi^{b_1}, \ldots, C^X(\psi_m)\varphi_{a_m} \psi^{b_m}). \]
    By the assumption $\mf{P}_{h,m}$ and the fact that $C^X$ has nonzero entries only in the top four rows, preserves the mod-$N$ degree, and the fact that $(\phi^j, C_{\ell}^X \varphi_a) \in \Q[X]_{\ell}$, we see that
    \[ \ab(\frac{Y}{t^N})^{h-1}\tilde{f}_{h,(\ba,\bb)}^{[0]} \in \Q[X]_{3h-3+m-\ab|\bb|}. \]

    In the graph sum formula for the action of $R^X$ and $\Omega^{\bA}$, we see there is a graph with a single genus $h$ vertex with $m$ insertions (corresponding to the largest stratum of $\ol{\msc{M}}_{h,m}$). The contribution of this graph is $f_{h,(\ba,\bb)}^{\bA}$. The contribution of any other graph will have the form
    \[ \bigotimes_{v} f^{\bA}_{g_v, n_v} \ab(\bigotimes_{i=1}^m \varphi_{a_i} \psi^{b_i} \otimes \bigotimes_{e} V_X(e)), \]
    where $V_X(e)$ was defined in~\Cref{cor:edgeX}. 

    We will now induct on $(h,m)$. In the base case $(h,m) = (0,3)$, the leading graph is the only graph, so the result follows directly (note here the different normalization conventions for $f^{\bA}$ and $f^{[0]}$). Now we assume the result for all $(h',m') \prec (h,m)$. Using the graph sum, we now count the degrees of all of the contributions.
    \begin{itemize}
        \item The total exponent of $\frac{Y}{t^N}$ is $\sum_{v} (g_v-1) + \ab|E| = h-1$ using~\Cref{cor:edgeX} and distributing the factors of $Y$ as in the proof of~\cite[Theorem 6.1]{polynomiality};
        \item The total degree in $X$ is at most
            \begin{align*}
                &\sum_v \ab(3g_v-3+n_v-\sum_{e\in E_v} b_{(e,v)} - \sum_{i \in L_v} b_i) + \sum_e (b_{ (e,v_1) }+b_{(e,v_2)}+1) \\
                ={}& 3 \ab(\sum_v g_v - \ab|V| + 3\ab|E|) + m - \sum_i b_i \\
                ={}& 3h-3+m - \ab|\bb|. \qedhere
            \end{align*}
    \end{itemize}
\end{proof}

\begin{lem}\label{lem:pimpliesq}
    If $\mf{P}_{h,m}$ holds for all $(h,m) \prec (g,n)$, then $\mf{Q}_{h,m}$ also holds for all $(h,s) \prec (g,n)$.
\end{lem}

\begin{proof}
    Recall that $\Omega^{[0]} = R^X.\Omega^{\bA}$. This implies that $f^{[0]}_{h,(\ba,\bb),(\ba',\bb')}$ can be written as a graph sum, where the contribution of a stable graph $\Gamma$ is given by the following:
    \begin{itemize}
        \item For each ordinary leg $\ell$, we insert $R^X(-\psi)^* \phi_a \psi^b$. By~\Cref{lem:matrixelementsRX}, this in fact becomes
            \[ (\varphi^{\ol{a} - m_{\ell}}, (-1)^{m_{\ell}}\psi^{b+m_{\ell}}(R^X_{m_{\ell}})^* \phi_a) \in \psi^{b+m_{\ell}} \ab(\frac{Y}{t^N})^{-\floor*{\frac{a}{N}}} \Q[X]_{m_{\ell}+\floor*{\frac{a}{N}}}. \]
            for a unique $m_{\ell}$ (here, note that the ancestor degree must be at most $3g_v-3+n_v < N$).
        \item For each special leg $\ell'$, we insert $R^X(-\psi)^* E_{a',b'}(\psi)$. Using~\Cref{lem:vanishingcoeffedge} and~\Cref{lem:degreecoeffedge}, this becomes
            \[ \psi^{b''} [z^{b''}] (\varphi^{a''}, R^X(-z)^* E_{a',b'}(z)) \in \psi^{b''} \ab(\frac{Y}{t^N})^{c_{E_{\ell'}}} \Q[X]_{b''+b'+1} \]
            for a unique $a'', b''$.
        \item At every edge, we insert the bivector $V_X$.
    \end{itemize}

    We now consider the total degree of the contributions from a graph $\Gamma$.
    \begin{itemize}
        \item The total exponent of $\frac{Y}{t^N}$ is given by
            \[ \sum_{\ell} -\floor*{\frac{a_{\ell}}{N}} + \sum_{\ell'} c_{E_{\ell'}} + \sum_v (1-g_v) +\sum_e (-1). \]
            Using the fact that
            \[ \sum_{\ell} (\bar{a}_{\ell} + b_{\ell}) + \sum_{\ell'} (a_{\ell'}'' + b_{\ell'}'') = m+n \]
            by~\Cref{cor:edgeX}, we obtain
            \begin{align*}
                \sum_{\ell'} c_{E_{\ell'}} + c + n - \floor*{\frac{\ba}{N}} ={}& \sum_{\ell'} \frac{a'_{\ell'} + b'_{\ell'} + a''_{\ell'} + b''_{\ell'}-N-2}{N} \\ 
                &+ \frac{\ab|\bar{\ba}| + \ab|\bb| - \ab|\ba'| - \ab|\bb'| -m+n}{N} + n \\
                ={}& \frac{\sum_{\ell'}(a''_{\ell'}+b''_{\ell'})+\ab|\bar{\ba}| + \bar|\bb| -nN - 2n -m+n }{N} \\
                ={}& 0,
            \end{align*}
            which implies that the total exponent is $-(c+n+h-1)$.
        \item The total degree in $X$ is at most 
            \begin{align*}
                &\sum_v \bigg(3g_v-3+n_v - \sum_{\ell}(b_{\ell}+m_{\ell})-\sum_{\ell'} b''_{\ell'} - \sum_e m_{(e,v)} \\
                &+ \sum_{\ell} \ab(m_{\ell}+\floor*{\frac{a_{\ell}}{N}}) + \sum_{\ell'} (b''_{\ell'}+b'_{\ell'}+1) \bigg) + \sum_e (m_{(e,v_1)}+ m_{(e,v_2)}+1) \\
                ={}& \ab(\sum_v 3g_v-3+n_v) - \ab|\bb| + \ab|\bb'| + \ab|E| + n + \floor*{\frac{\ba}{N}} \\
                ={}& 3h-3+m+2n+\floor*{\frac{\ba}{N}} - \ab|\bb| + \ab|\bb'|. \qedhere
            \end{align*}
    \end{itemize}
\end{proof}

\begin{thm}
    For all $\ba, \bb \in \{0,\ldots,N+3\}^n$, we have
    \[ \ab(\frac{Y}{t^N})^{g-1+c} f_{g,(\ba,\bb)}^{[0]} \in \Q[X]_{3g-3+n+\floor*{\frac{\ba}{N}} - \ab|\bb|}. \]
\end{thm}

\begin{proof}
    We will induct on $(g,n)$ under the ordering $\prec$. The base case is $(g,n) = (0,3)$. In this case, there is only one graph with a single vertex. Because $\dim \ol{\msc{M}}_{g,n} = 0$, no ancestor insertions are allowed, and so we calculate
    \begin{align*}
        \ab(\frac{Y}{t^N})^{0-1+c} f_{0,(\ba,\mbf{0})}^{[0]} &= \ab(\frac{Y}{t^N})^c I_0^2 I_{11}^2 I_{22} \frac{t^N}{Y} \cdot \mr{const} \cdot q^c \\
        &= \mr{const} \cdot X^c.
    \end{align*}
    Here, we use the fact that $c = \floor*{\frac{\ba}{N}}$ and the computation of genus-zero three-point functions in~\cite[\S 2.4]{polynomiality}.

    We now assume the desired polynomiality result for $(h,m) \prec (g,n)$. This implies $\mf{P}_{h,m}$ and thus $\mf{Q}_{h,m}$ for all $(h,m) \prec (g,n)$ by~\Cref{lem:pimpliesq}. We may also assume that $c = \floor*{\frac{\ba}{N}}$ by~\Cref{cor:criterionforc}. 

    We will now consider the $[0,1]$ theory. By~\cite[Theorem 4.1]{polynomiality}, $f_{g,(\ba,\bb)}^{[0,1]}$ is a polynomial in $q$ of degree at most $g-1+\frac{3g-3 + \ab|\ba|}{N}$. By~\Cref{lem:vanishing0nospecial}, this becomes
    \begin{align*}
        g-1+\frac{3g-3+\ab|\ba|+n-\ab|\bar{\ba}|-\ab|\bb|}{N} &= g-1+\frac{3g-3+\ab|\ba|+n-\ab|\ba|-\ab|\bb|+N\floor*{\frac{\ba}{N}}}{N} \\
        &= g-1+c + \frac{3g-3+n-\ab|\bb|}{N}.
    \end{align*}
    Because $0\leq \ab|\bb| \leq 3g-3+n$ and $N \gg 3g-3+n$, we see that the degree is in fact at most $g-1+c$. Multiplying by $\ab(\frac{Y}{t^N})^{g-1+c}$, we see that
    \[ \ab(\frac{Y}{t^N})^{g-1+c} f_{g,(\ba,\bb)}^{[0,1]} \in \Q[X]_{g-1+c}. \]
    By~\Cref{cor:criterionforc}, this satisfies the desired degree bound.

    We will now apply the bipartite graph decomposition from~\Cref{thm:01bipartite}. There is a leading bipartite graph with only a single level $0$ vertex. We need to prove the degree estimate for the non-leading graphs. Applying $\mf{Q}_{h,m}$ at level $0$ and~\Cref{lem:polynomiality1theory} at level $1$, we now count the total degree contribution of a bipartite graph $\Lambda$.
    \begin{itemize}
        \item The total exponent of $\frac{Y}{t^N}$ is
            \begin{align*}
                & g-1+c - \sum_{v \in V_0} (g_v-1+c_v+\ab|L'_v|) - \sum_{v \in V_1} (g_v-1+c_v) \\
                ={}& g-1 + \frac{\ab|\ba|+\ab|\bb|-n}{N}\\ 
                &- \sum_{v \in V_0} \ab(g_v-1+\frac{\ab|\ba_v| + \ab|\bb|_v -\ab|\ba'_v| - \ab|\bb'_v| - \ab|L_v| + \ab|E_v|}{N} + \ab|E_v|) \\
                &- \sum_{v \in V_1} \ab(g_v-1+\frac{\ab|\ba_v| + \ab|\bb|_v -\ab|\ba'_v| - \ab|\bb'_v| - \ab|L_v| - \ab|E_v|}{N}) \\
                ={}& g-1 + \ab|E| - \sum_{v \in V} (g_v-1) \\
                ={}& 0.
            \end{align*}
        \item The total degree in $X$ is
            \begin{align*}
                \sum_{v \in V_0} \ab(3g_v-3+n_v + \ab|E_v|+\floor*{\frac{\ba_v}{N}} - \ab|\bb_v| + \ab|\bb'_v|) \\
                + \sum_{v \in V_1} \ab(3g_v - 3 + n_v + \floor*{\frac{\ba_v}{N}} - \ab|\bb_v| - \ab|\bb'_v|)
            \end{align*}
            Because a stable graph describes a codimension $\ab|E|$ stratum in $\ol{\msc{M}}_{g,n}$, this is at most $3g-3+n+\floor*{\frac{\ba}{N}} - \ab|\bb|$. \qedhere
    \end{itemize}
\end{proof}

Applying~\Cref{lem:Apolyassuming0poly}, we obtain the following.
\begin{cor}\label{cor:amodelfeynman}
    For all $(g,n)$ such that $2g-2+n > 0$, all $\ba \in \{0,1,2,3\}^n$, and all $\bb \in \Z_{\geq 0}^n$, we have
    \[ f_{g,(\ba,\bb)}^{\bA} \in \Q[X]_{3g-3+n-\ab|\bb|}. \]
\end{cor}

\subsection{Choice of gauge}%
\label{sub:Choice of gauge}

Note that
\[ R^{\bA,\G}(z)^{-1} = R^{\bA}(z)^{-1} \G(z)^{-1}, \]
where we compute
\[ \G(z)^{-1} = I - \begin{pmatrix}
    0 & zc_{11} & z^2c_2 & -z^3 (c_{11} c_2 + c_3) \\
    & 0 & z c_{12} & -z^2 (c_{11} c_{12} + c_2) \\
    & & 0 & c_{11} \\
    & & & 0
\end{pmatrix}.
\]
Note that this is a symplectic matrix, so we have an equality
\[ \Omega^{\bA, \G} = \G.\Omega^{\bA} \]
of CohFTs.

\begin{thm}
    For all $(g,n)$ such that $2g-2+n > 0$, all $\ba \in \{0,1,2,3\}^n$, and all $\bb \in \Z_{\geq 0}^n$, we have
    \[ f_{g,(\ba,\bb)}^{\bA,\G} \in \Q[X]_{3g-3+n-\ab|\bb|}. \]
\end{thm}

\begin{proof}
    We will write the stable graph sum formula for $\Omega^{\bA, \G}$ as the $\G$-action on $\Omega^{\bA}$. The contribution of a stable graph $\Gamma$ is given by the following assignments:
    \begin{itemize}
        \item At each leg, we place $\G(-z)^* \varphi_a \psi^b = \sum_m \G_m^* (-\psi)^m \varphi_a \psi^b$;
        \item At each edge, we place
            \begin{align*}
                V^{\G} &\coloneqq \sum_{i=1}^3 \frac{\varphi_i \otimes \varphi^i - \G(-\psi)^* \varphi_i \otimes \G(-\psi')^* \varphi^i}{\psi+\psi'} \\
                &= Y^{-1} \sum_{a,b} V_{ab}^{\G} \psi^a (\psi')^b.
            \end{align*}
    \end{itemize}
    Here, note that $\G_m^*$ has degree $m$ in $X$ and $V_{ab}$ has degree $a+b+a$ in $X$ by the assumption on the gauge in~\Cref{defn:propogators}.

    We now compute the total degree of the contribution. The total exponent of $\frac{Y}{t^N}$ in the contribution of $\Gamma$ to $\Omega_{g,n}^{\bA, \G}$ is
    \[ -\ab|E| - \sum_{v \in V} (g_v-1) = -(g-1). \]
    On the other hand, the total degree in $X$ is at most
    \begin{align*}
        &\sum_v \ab(3g_v-3+n_v - \sum_{\ell \in L_v} (m_{\ell} + b_{\ell}) - \sum_{e \in E_v} m_{(e,v)})\\
        &+ \sum_e (m_{(e,v_1)} + m_{(e,v_2)}+1) + \sum_{\ell} m_{\ell} \\
        ={}& 3g-3+n - \ab|\bb|. \qedhere
    \end{align*}
\end{proof}

\section{The B-model Feynman rule}%
\label{sec:The B-model Feynman rule}

In this section, we will define the B-model Feynman rule and prove that it equals the A-model Feynman rule. From now on, we will make the specialization $t^N = -1$. This makes $q' = q$ and makes
\[ f_{g,(\ba,\bb)}^{\bA,\G} = (p_kY)^{g-1} \int_{\ol{\msc{M}}_{g,n}} \Omega_{g,n}^{\bA,\G}(\varphi_{a_1} \psi_1^{b_1}, \ldots, \varphi_{a_n} \psi_n^{b_n}). \]

\begin{conv}
    In this section, we will omit all superscripts of $\G$, so $E_{**}$ stands for $E_{**}^{\G}$, $R^{\bA}$ stands for $R^{\bA, \G}$, and so on.
\end{conv}

\subsection{B-model geometric quantization}%
\label{sub:B model geometric quantization}

We will first express the physics Feynman rule using geometric quantization. Note that the Givental formalism is a type of geometric quantization, so we will be able to compare the A-model and B-model quantizations.

Consider the vector space 
\[\msc{H}_S \coloneqq T^*(zH^0(Z) \oplus H^2(Z)) = \on{span}\{  -\varphi_2 z^{-1}, \varphi_3 z^{-2},\varphi_1, \varphi_0 z \} \] with the symplectic form
\[ \frac{1}{p_k Y} \Res_{z=0} (f(-z),g(z)) = \begin{pmatrix}
    0 & I \\
    -I & 0
\end{pmatrix}
. \]
Then define
\begin{align*}
    R^{\bB} &\coloneqq R^{\bA}|_{\msc{H}_S} \\
    &= \begin{pmatrix}
        1 & -E_{\psi} \\
        0 & 1 \\
        -E_{\varphi\varphi} & -E_{\varphi\psi} & 1 \\
        E_{1\varphi\psi} & E_{1\psi^2} & E_{\psi} & 1
    \end{pmatrix} \\
    &\eqqcolon \begin{pmatrix}
        A & B \\
        C & D
    \end{pmatrix}.
\end{align*}
Because $\msc{H}_S$ is a symplectic vector space, we will write elements as
\[ (\mbf{p}, \mbf{x}) \coloneqq  - p_x \varphi_2 z^{-1} +p_y \varphi_3 z^{-2} + x \varphi_1 + y \varphi_0 z. \]

\begin{defn}
Following~\cite{geomquantgwtheory}, we define the geometric quantization $\wh{R}^{\bB}$ by the Gaussian integral
\[ \wh{R}^{\bB}F(\hbar, \mbf{x}) \coloneqq \log \int_{\R^4} e^{\frac{1}{\hbar} ( \mbf{Q}(\bx', \bp') - \bx' \cdot \bp' ) + F(\hbar, \bx')} \d{\bx'} \d{\bp'}. \]
Here, $\mbf{Q}(\bx', \bp')$ is defined by the formula
\begin{align*}
    \mbf{Q}(\bx', \bp') &\coloneqq (D^{-1}\bx') \cdot \bp' - \frac{1}{2} (D^{-1}C \bp') \cdot \bp' \\
    &= (\bp')^T \begin{pmatrix}
        1 & \\
        -E_{\psi} & 1
    \end{pmatrix} \bx' + \frac{1}{2} (\bp')^T \begin{pmatrix}
        E_{\varphi\varphi} & E_{\varphi\psi} \\
        E_{\varphi\psi} & E_{\psi\psi}
    \end{pmatrix}\bp'
\end{align*}
and $(\bp', \bx')$ are coordinates on $\R^4 = \R^2 \times \R^2$.
\end{defn}

Following~\cite[\S 3.4]{geomquantgwtheory}, a standard argument involving the Fourier transform gives us the operator form
\[ \wh{R}^{\bB} F(\hbar, \bx) = \log \ab(e^{-\frac{\hbar}{2} \begin{psmallmatrix}
    \partial_x & \partial_y
\end{psmallmatrix}
DC^T \begin{psmallmatrix}
    \partial_x \\ \partial_y
\end{psmallmatrix}
} e^{F(\hbar, D^{-1}\bx)}). \]
Now define $\tilde{E}_{\varphi\varphi}$, $\tilde{E}_{\varphi\psi}$, and $\tilde{E}_{\psi\psi}$ by
\[ -DC^T \coloneqq \begin{pmatrix}
    \tilde{E}_{\varphi\varphi} & \tilde{E}_{\varphi\psi} \\
    \tilde{E}_{\varphi\psi} & \tilde{E}_{\psi\psi}
\end{pmatrix}.
\]
Then if we define
\[ V^{\bB}(\partial_{\bx}, \partial_{\bx}) \coloneqq \frac{1}{2} \tilde{E}_{\varphi\varphi} \pdv[order=2]{}{x} + \tilde{E}_{\varphi\psi} \pdv{}{x,y} + \frac{1}{2} \tilde{E}_{\psi\psi} \pdv[order=2]{}{y}, \]
the operator form of the quantization action becomes
\begin{equation}\label{eqn:operatorbmodel}
    \wh{R}^{\bB} F(\hbar, \bx) = \log \ab(e^{\hbar V^{\bB}(\partial_{\bx}, \partial_{\bx})} e^{F(\hbar, D^{-1}\bx)}). 
\end{equation}

\begin{defn}
    Define the (normalized) Gromov-Witten correlator of $Z$ by the formula
    \[ P_{g,m,n} \coloneqq \frac{(p_k Y)^{g-1}}{I_0^{2g-2+m+n}} \ab< \varphi_1^{\otimes m}, (\varphi_0 \psi)^{\otimes n}>_{g,m+n}^Z \]
    when $(g,m) \neq (1,0)$ and define
    \[ P_{1,0,n} \coloneqq (n-1)! \ab(\frac{\chi(Z)}{24} - 1). \]
\end{defn}
 
\begin{defn}
    Define the master B-model Gromov-Witten potential function by
    \[ P^{\bB}(\hbar, x, y) \coloneqq \sum_{g,m,n} \hbar^{g-1} \frac{x^m y^n}{m!n!} P_{g,m,n}. \]
    Then, define the master B-model potential function by
    \[ f^{\bB}(\hbar, x, y) \coloneqq \wh{R}^{\bB} P^{\bB}(\hbar, x, y) \eqqcolon \sum_{g,m,n} \hbar^{g-1} f_{g,m,n}^{\bB}. \]
\end{defn}

By~\cite[Theorem 10]{geomquantgwtheory}, we can also compute $f^{\bB}$ by the construction in~\Cref{defn:bmodelfeynman}.

\subsection{Factorization of the quantization action}%
\label{sub:Factorization of the quantization action}

We will factor the quantization action~\Cref{eqn:operatorbmodel} into the change of variables and the application of differential operators. Observe that $D^{-1}\bx = (x,y-E_{\psi} x)$. Then the transformation
\[ F(\hbar, x,y) \mapsto F(\hbar, x, y-E_{\psi} x) \]
is given by quantizing the matrix
\[ \msc{E}^{\bB} \coloneqq \begin{pmatrix}
    1 & -E_{\psi} \\
     & 1 \\
    & & 1 \\
    & & E_{\psi} & 1
\end{pmatrix}.
\]
Then we compute
\begin{align*}
    \tilde{P}^{\bB} &\coloneqq \wh{\msc{E}}^{\bB} P^{\bB}(\hbar, x, y) \\
    &= P^{\bB}(\hbar, x, y-E_{\psi}x) \\
    &= \sum_{g,m,n} \hbar^{g-1} \frac{x^m y^n}{m!n!} \tilde{P}_{g,m,n},
\end{align*}
where
\[ \tilde{P}_{g,m,n} \coloneqq \frac{(p_k Y)^{g-1}}{I_0^{2g-2+m+n}}\ab<(\varphi_1 - E_{\psi} \varphi_0\psi)^{\otimes m}, (\varphi_0 \psi)^{\otimes n}>_{g,m+n}^{Z}. \]
We then compute
\[ \tilde{R}^{\bB} \coloneqq R^{\bB} (\msc{E}^{\bB})^{-1} = \begin{pmatrix}
    1 & 0 \\
    0 & 1 \\
    -\tilde{E}_{\varphi\varphi} & -\tilde{E}_{\varphi\psi} & 1 \\
    -\tilde{E}_{\varphi\psi} & -\tilde{E}_{\psi\psi} & 0 & 1
\end{pmatrix},
\]
so we see that $f^{\bB}(\hbar, x, y) = \wh{\tilde{R}^{\bB}} \tilde{P}^{\bB}(\hbar, x, y)$.

\subsection{Modification of the A-model quantization}%
\label{sub:Modification of the A-model quantization}

Note that in the graph sum formula for $\Omega^{\bA}$, the edge contribution is given by
\begin{align*}
    V_{\bA} 
    ={}& V_{\bB} 
    + E_{\psi} (\varphi_0 \otimes \varphi_2 + \varphi_2 \otimes \varphi_0) + E_{1\varphi\psi} (\varphi_0 \otimes \varphi_1 \psi' + \varphi_1 \psi \otimes \varphi_0) \\ 
    &+ E_{1\psi^2} (\varphi_0 \otimes \varphi_0 (\psi')^2 + \varphi_0 \psi^2 \otimes \varphi_0).
\end{align*}
In order to prove that the A-model and B-model Feynman rules are equivalent, we need to analyze the contributions of the three extra terms. We will begin with the terms $E_{\psi} (\varphi_0 \otimes \varphi_2 + \varphi_2 \otimes \varphi_0)$ and study a parallel construction to the modified B-model quantization.

The parallel construction in the A-model is to consider the matrix\footnote{This should also implicitly have a $\Psi^{-1}$ on the right, but we will omit it.}
\[ \msc{E}^{\bA} \coloneqq I + z \begin{pmatrix}
    0 & E_{\psi} \\
    & 0 \\
    & & 0 & E_{\psi} \\
    & & & 0
\end{pmatrix}
\]
and the factorization
\[ \tilde{R}^{\bA}(z) \coloneqq R^{\bA}(z) \msc{E}^{\bA}(z)^{-1}. \]
Recall that matrices are written in the basis $\varphi_0, \varphi_1, \varphi_2, \varphi_3$.

\begin{defn}
    Define the CohFT
    \[ \tilde{\Omega}^Z \coloneqq \msc{E}^{\bA} \Omega^{Z} \]
    and the potential
    \[ \tilde{P}^{\bA}(\hbar, \bt) \coloneqq \sum_{g,n} \frac{\hbar^{g-1}}{n!} (p_k Y)^{g-1} \int_{\ol{\msc{M}}_{g,n}} \tilde{\Omega}^Z_{g,n}(\bt^{\otimes n}) \]
    for the coordinate
    \[ \bt = x \varphi_1 + y \varphi_0 \psi + a \varphi_1 \psi + b \varphi_0 \psi^2 + c \varphi_0. \]
    Specializing to the case $\bt = (\bx, 0)$, we define
    \[ \tilde{P}^{\bA}(\hbar, x, y) \coloneqq \tilde{P}^{\bA}(\hbar, x\varphi_0 + y \varphi_0 \psi). \]
\end{defn}

\begin{lem}
    We have the identity
    \[ \tilde{P}^{\bA}(\hbar, x, y) = \tilde{P}^{\bB}(\hbar, x, y) - \log(1-y). \]
\end{lem}

\begin{proof}
    By~\cite[Theorem 5]{invtauteqns}, the $R$-matrix action preserves all tautological equations on the moduli spaces of curves, in particular the string and dilaton equations, so we can apply the dilaton equation to $\tilde{\Omega}^Z$. We also note that changing $\tilde{P}^{\bB}(\hbar, x, y)$ to $\tilde{P}^{\bB}(\hbar, x, y)$ is the same as replacing $P_{1,0,n}$ by the actual Gromov-Witten invariant
    \[ \frac{1}{I_0^n} \ab<(\varphi_0\psi)^{\otimes n}>_{1,n}^Z = (n-1)! \frac{\chi(Z)}{24}, \]
    so applying the dilaton equation to both sides (where $\tilde{P}$ denotes either $\tilde{P}^{\bA}$ or $\tilde{P}^{\bB}$) yields
    \begin{align*}
        \pdv{}{y} \tilde{P}(\hbar, x, y) &= \pdv{}{y}\ab( \sum_{(g,m,n) \neq (1,0,1)} \hbar^{g-1} \frac{x^m y^n}{m!n!} \tilde{P}_{g,m,n} + \frac{\chi(Z)}{24} y ) \\
        &= \sum_{g,m,n} \hbar^{g-1} \frac{x^m y^{n-1}}{m!(n-1)!} \tilde{P}_{g,m,n} + \frac{\chi(Z)}{24} \\
        &= \sum_{g,m,n} \hbar^{g-1} \frac{x^m y^{n-1}}{m!(n-1)!} (2g-2+m+n-1) \tilde{P}_{g,m,n-1} + \frac{\chi(Z)}{24} \\
        &= \ab(2\hbar \pdv{}{\hbar} + x \pdv{}{x} + y \pdv{}{y}) \tilde{P} + \frac{\chi(Z)}{24}.
    \end{align*}
    Therefore, we only need to prove that $\tilde{P}^{\bA}(\hbar, x, 0) = \tilde{P}^{\bB}(\hbar, x, 0)$.

    We will now consider the graph sum formula for the definition of $\tilde{\Omega}^Z$. The edge contributions are given by $E_{\psi}(\varphi_0 \otimes \varphi_2 + \varphi_2 \otimes \varphi_0)$. All insertions at legs are $\msc{E}^{\bA}(-\psi)^* \varphi_1 = \varphi_1 - E_{\psi} \varphi_0 \psi$, and each vertex contributes a Gromov-Witten correlator of $Z$. Applying the string equation, dilaton equation, divisor equation, and virtual dimension constraints, if the stable graph $\Gamma$ has at least one edge and one vertex of $g > 0$, then its contribution vanishes. By a similar argument, any vertex with more than two edges has vanishing contribution. Therefore, we have a decomposition
    \[ \tilde{P}^{\bA}(\hbar, x, 0) = P^{\bA}(\hbar, x, -E_{\psi}x) + P_1^{\sloop}(x), \]
    where the first term comes from the leading graphs with a single genus $g$ vertex and the second comes from loops of genus zero vertices, which contribute to the genus $1$ potential. Every vertex must have at least one $\varphi_1$ insertion by dimension reasons, and the string and dilaton equations imply that there is exactly one $\varphi_1$ insertion. Therefore, we use the dilaton equation to remove the $-E_{\psi}\varphi_0 \psi$ insertions and compute
    \begin{align*}
        P_1^{\sloop}(x) &= \sum_{\substack{\Gamma \text{ loop} \\ \ab|V| = m \\ \ab|L| = m+n}} \frac{x^{m+n}}{m!n!}\frac{\mr{Cont}_{\Gamma}}{\ab|\Aut \Gamma|} \\
        &= \sum_{m > 0} \frac{(m-1)!}{m!} (E_{\psi}x)^m \prod_{i=1}^m \sum_{n_i \geq 0} (-E_{\psi}x)^{n_i} \\
        &= \sum_{m > 0} \frac{1}{m} \ab(\frac{E_{\psi}x}{1+E_{\psi}x})^m \\
        &= -\log \ab(1- \frac{E_{\psi}x}{1+E_{\psi}x}) \\
        &= \log (1+E_{\psi}x).
    \end{align*}
    Here, we have used the following combinatorial facts:
    \begin{itemize}
        \item There are $(m-1)!$ ways to arrange $m$ vertices in a loop;
        \item For any partition $n = n_1 + \cdots + n_m$ of the $\varphi_0 \psi$ insertions, the number of possible assignments is $\frac{n!}{n_1! \cdots n_m!}$;
        \item Applying the dilaton equation to any vertex with $n_i$ insertions of $\varphi_0 \psi$ and one insertion of $\varphi_1$ produces a factor of $n!$.
    \end{itemize}
    We conclude that 
    \begin{align*}
        \tilde{P}^{\bA}(\hbar, x, 0) &= P^{\bA}(\hbar, x, -E_{\psi}x) + \log(1+E_{\psi}x) \\
        &= P^{\bB}(\hbar, x, -E_{\psi}x) \\
        &= \tilde{P}^{\bB}(\hbar, x, 0). \qedhere
    \end{align*}
\end{proof}

\subsection{Equality of A-model and B-model potentials}%
\label{sub:Equality of A-model and B-model potentials}

A direct computation yields
\[ \tilde{R}^{\bA}(z)^{-1} = I - \begin{pmatrix}
    0 & 0 & z^2 \tilde{E}_{\varphi\psi} & z^3 \tilde{E}_{\psi\psi} \\
    & 0 & z \tilde{E}_{\varphi\varphi} & z^2 \tilde{E}_{\varphi\psi} \\
    & & 0 & 0 \\
    & & & 0
\end{pmatrix}.
\]
By~\cite[Proposition 7.3]{virasorofanotoric}, we see that
\[ e^{f^{\bA}(\hbar, x, y)} = e^{\hbar V^{\bA}(\partial_{\bt}, \partial_{\bt})} e^{\tilde{P}^{\bA}(\hbar, x, y)}, \]
where
\begin{align*}
    V^{\bA}(\partial_{\bt}, \partial_{\bt}) &\coloneqq V^{\bB}(\partial_{\bt}, \partial_{\bt}) - \tilde{E}_{\varphi\psi} \pdv{}{a,c} - \tilde{E}_{\psi\psi} \pdv{}{b,c} \\
    &\eqqcolon V^{\bB}(\partial_{\bt}, \partial_{\bt}) + V^{\extra}(\partial_{\bt}, \partial_{\bt}). 
\end{align*}

We will first prove some technical lemmas about the geometric quantization formalism. This is

\begin{lem}
    We have the equality
    \[ e^{\tilde{P}^{\bA}(\hbar, x, y, a, b, c)} = e^{\frac{c}{1-y}\ab(a \pdv{}{x} + b \pdv{}{y})} e^{\tilde{P}^{\bA}(\hbar, x, y)}. \]
\end{lem}

\begin{proof}
    The operator form of the string equation\footnote{If we write $\varphi(z) = \sum_{i=0}^3 t_i^j \varphi_i z^j$, the operator form of the string equation (see~\cite{symplfrob}) is
    \[ \pdv{}{t_0^0} \msc{D} = \frac{1}{2} (\bt^0, \bt^0) + \sum_{j=0}^{\infty} \sum_{i=0}^3 t_i^{j+1} \pdv{}{t_i^j} \msc{D}. \]
    Because there are no $\varphi_2$ or $\varphi_3$ insertions, the unstable term disappears. We then set $t_1^0 = x$, $t_1^1 = a$, $t_0^0 = c$, $t_0^1 = y$, and $t_0^2 = b$. Here, $\msc{D} = e^{\tilde{P}^{\bA}}$.} is
    \[ \pdv{}{c} e^{\tilde{P}^{\bA}(\hbar, x, y, a, b, c)} = \ab(a \pdv{}{x} + b \pdv{}{y} + y\pdv{}{v})e^{\tilde{P}^{\bA}(\hbar, x, y, a, b, c)}. \]
    By virtual dimension reasons, we obtain the initial condition
    \[ \tilde{P}^{\bA}(\hbar, x, y, a, b, 0) = \tilde{P}^{\bA}(\hbar, x, y). \]
    The desired result follows from the computation
    \begin{align*}
        (1-y)\pdv{}{c}&\ab( e^{\frac{c}{1-y}\ab(a \pdv{}{x} + b \pdv{}{y})} e^{\tilde{P}^{\bA}(\hbar, x, y)} ) \\
        &= (1-y)\pdv{}{c} \sum_{n=0}^{\infty} \frac{1}{n!} \ab(\frac{c}{1-y})^n \ab(a \pdv{}{x} + b \pdv{}{y})^n e^{\tilde{P}(\hbar, x, y)} \\
        &= \sum_{n=0}^{\infty} \frac{1}{(n-1)!} \ab(\frac{c}{1-y})^{n-1} \ab(a \pdv{}{x} + b \pdv{}{y})^n e^{\tilde{P}^{\bA}(\hbar, x, y)} \\
        &= \ab(a \pdv{}{x} + b\pdv{}{y}) \ab(e^{\frac{c}{1-y}\ab(a \pdv{}{x} + b \pdv{}{y})} e^{\tilde{P}^{\bA}(\hbar, x, y)}). \qedhere
    \end{align*}
\end{proof}

\begin{thm}\label{thm:allgenusmirror}
    We have the equality
    \[ f^{\bA}(\hbar, x, y) = f^{\bB}(\hbar, x, y) - \log(1-y). \]
\end{thm}

\begin{proof}
    Our goal is to compute the function
    \begin{align*}
        e^{f^{\bA}(\hbar, x, y)} ={}&
        e^{\hbar(V^{\bB}(\partial_{\bt}, \partial_{\bt}) + V^{\extra}(\partial_{\bt}, \partial_{\bt}))} e^{\tilde{P}^{\bA}(\hbar, x, y, a, b, c)}\Big|_{a,b,c=0} \\
        ={}& e^{\hbar(V^{\bB}(\partial_{\bt}, \partial_{\bt})+V^{\extra}(\partial_{\bt}, \partial_{\bt}))} e^{\frac{c}{1-y}\ab(a\pdv{}{x}+b\pdv{}{y})}e^{\tilde{P}^{\bA}(\hbar, x, y)} \Big|_{a,b,c=0} \\
        ={}& e^{\hbar(V^{\bB}(\partial_{\bt}, \partial_{\bt})+V^{\extra}(\partial_{\bt}, \partial_{\bt}))} e^{\frac{c}{1-y}\ab(a\pdv{}{x}+b\pdv{}{y})}\frac{e^{\tilde{P}^{\bB}(\hbar, x, y)}}{1-y} \Bigg|_{a,b,c=0}. 
    \end{align*}
    We will first consider the contribution of $\frac{c}{1-y} \ab(a\pdv{}{x}+b\pdv{}{y})$ and $V^{\extra}(\partial_{\bt}, \partial_{\bt})$. For any function $\msc{D}(x,y)$, we compute
    \begin{align*}
        & e^{\hbar V^{\extra}(\partial_{\bt}, \partial_{\bt})} e^{\frac{c}{1-y}\ab(a\pdv{}{x}+b\pdv{}{y})} \msc{D}(x,y) \\
        ={}& \sum_{m,n} \frac{(-\hbar)^m}{m!n!}  \ab(\tilde{E}_{\varphi\psi}\pdv{}{a,c} + \tilde{E}_{\psi\psi}\pdv{}{b,c})^m \ab(\frac{c}{1-y} \ab(a\pdv{}{x}+b\pdv{}{y}))^n \msc{D}(x,y) \Big|_{a,b,c=0} \\
        ={}& \sum_{n} \frac{(-\hbar)^n}{(n!)^2}  \ab(\tilde{E}_{\varphi\psi}\pdv{}{a,c} + \tilde{E}_{\psi\psi}\pdv{}{b,c})^n \ab(\frac{c}{1-y} \ab(a\pdv{}{x}+b\pdv{}{y}))^n \msc{D}(x,y) \Big|_{a,b,c=0} \\
        ={}& \sum_n \frac{(-\hbar)^n}{n!}\ab(\tilde{E}_{\varphi\psi}\pdv{}{a} + \tilde{E}_{\psi\psi}\pdv{}{b})^n \ab(\frac{1}{1-y} \ab(a\pdv{}{x}+b\pdv{}{y}))^n \msc{D}(x,y) \Big|_{a,b,c=0} \\
        ={}& \sum_n \ab(\frac{(-\hbar)^n}{1-y} \ab(\tilde{E}_{\varphi\psi}\pdv{}{x} + \tilde{E}_{\psi\psi} \pdv{}{y}))^n \msc{D}(x,y).
    \end{align*}
    Now set $E^{\extra}(\partial_{\bt}) \coloneqq \frac{-\hbar}{1-y}\ab(\tilde{E}_{\varphi\psi}\pdv{}{x} + \tilde{E}_{\psi\psi}\pdv{}{y})$.

    To deal with the contribution of $V^{\bB}(\partial_{\bt}, \partial_{\bt})$, we compute
    \begin{align*}
        e^{-\hbar V^{\bB}(\partial_{\bt}, \partial_{\bt})}&(1-y)e^{\hbar V^{\bB}(\partial_{\bt}, \partial_{\bt})} \\ 
        &= \ab(\sum_{m\geq 0} \frac{(-\hbar V^{\bB}(\partial_{\bt}, \partial_{\bt}))^m}{m!})(1-y) \ab(\sum_{n\geq 0} \frac{(-\hbar V^{\bB}(\partial_{\bt}, \partial_{\bt}))^n}{n!}) \\
        &= \sum_{n \geq 0} \sum_{\ell+m = n} \frac{(-\hbar V^{\bB}(\partial_{\bt}, \partial_{\bt}))^{\ell}}{\ell!}(1-y) \frac{(\hbar V^{\bB}(\partial_{\bt}, \partial_{\bt}))^m}{m!} \\
        &= \sum_{n \geq 0}\frac{1}{n!} \ab(\ab[-,V^{\bB}(\partial_{\bt}, \partial_{\bt})])^n (1-y) \\
        &= (1-y) + \hbar\ab(\tilde{E}_{\varphi\psi} \pdv{}{x} + \tilde{E}_{\psi\psi} \pdv{}{y}).
    \end{align*}
    Dividing by $1-y$, we see that
    \[ (1-y)^{-1} e^{-\hbar V^{\bB}(\partial_{\bt}, \partial_{\bt})}(1-y)e^{\hbar V^{\bB}(\partial_{\bt}, \partial_{\bt})} = 1-E^{\extra}(\partial_{\bt}), \] 
    or in other words that
    \[  e^{-\hbar V^{\bB}(\partial_{\bt}, \partial_{\bt})}(1-y)^{-1}e^{\hbar V^{\bB}(\partial_{\bt}, \partial_{\bt})} (1-y) = \sum_{n \geq 0} E^{\extra}(\partial_{\bt})^n.\]

    Putting all of this together, we see that
    \begin{align*}
        e^{f^{\bA}(\hbar, x, y)}
        ={}& e^{\hbar(V^{\bB}(\partial_{\bt}, \partial_{\bt})+V^{\extra}(\partial_{\bt}, \partial_{\bt}))} e^{\frac{c}{1-y}\ab(a\pdv{}{x}+b\pdv{}{y})}\frac{e^{\tilde{P}^{\bB}(\hbar, x, y)}}{1-y} \Bigg|_{a,b,c=0} \\
        ={}& e^{\hbar V^{\bB}(\partial_{\bt}, \partial_{\bt})} \sum_{n \geq 0} E^{\extra}(\partial_{\bt})^n (1-y)^{-1} e^{\tilde{P}^{\bB}(\hbar, x, y)} \\
        ={}& (1-y)^{-1} e^{\hbar V^{\bB}(\partial_{\bt}, \partial_{\bt})} e^{\tilde{P}^{\bB}(\hbar, x, y)} \\
        ={}& \frac{e^{f^{\bB}(\hbar, x, y)}}{1-y}.
    \end{align*}
    The desired result follows by taking logarithms.
\end{proof}

\begin{cor}[B-model Feynman rule]\label{cor:bmodelfeynman}
    For any $g,m,n$, we have
    \[ f^{\bB}_{g,m,n} \in \Q[X]_{3g-3+m}. \]
\end{cor}

\begin{proof}
    By~\Cref{thm:allgenusmirror}, we see that $f^{\bA}_{g,m,n} = f^{\bB}_{g,m,n} + \delta_{g,1}\delta_{m,0} (n-1)!$. The result then follows from~\Cref{cor:amodelfeynman} by choosing $\ba = (1^m 0^n)$ and $\bb = (0^n 1^m)$.
\end{proof}

\section{Anomaly equations}%
\label{sec:Holomorphic anomaly equations}

Our goal is to prove the following theorem.
\begin{thm}\label{thm:hae}
    The $P_{g,m}$ satisfy the differential equations
    \begin{align}
    - \partial_{A} P_g = \frac{1}{2} \ab(P_{g-1,2} + \sum_{g_1+g_2 = g} P_{g_1, 1} P_{g_2, 2}), \label{eqn:hae1} \\
    \ab(-2 \partial_{A} + \partial_{B} + (A+2B) \partial_{B_2} - \ab((B-X)(A+2B)-B_2-r_0 X)\partial_{B_3})P_g = 0\label{eqn:hae2}.
    \end{align}
\end{thm}

In order to prove this, we introduce a smaller ring of modified generators which contains the $P_g$. 
\begin{defn}
Define the modified generators
\begin{align*}
    \msc{E}_1 \coloneqq \tilde{E}_{\varphi\varphi}, \qquad
    \msc{E}_2 \coloneqq \tilde{E}_{\varphi\psi}, \qquad
    \msc{E}_3 \coloneqq \tilde{E}_{\psi\psi} 
\end{align*}
and then set
\[ \tilde{\msc{R}} \coloneqq \Q[\msc{E}_1, \msc{E}_2, \msc{E}_3, X] \subset \msc{R}. \]
\end{defn}

\begin{rmk}
    The ring $\tilde{\msc{R}}$ is in fact invariant under the choice of gauge. Direct computation yields
    \begin{align*}
        \msc{E}_1^{\G} &= \msc{E}_1^{\bO} + c_{12}, \\
        \msc{E}_2^{\G} &= \msc{E}_2^{\bO} + c_{11} \msc{E}_1^{\bO} + c_{11}c_{12} + c_2, \\
        \msc{E}_3^{\G} &= \msc{E}_3^{\bO} + 2c_{11} \msc{E}_2^{\bO} + c_{11}^2 \msc{E}_1^{\bO} + c_{11}^2 c_{12} + 2 c_{11}c_2 + c_3,
    \end{align*}
    so we will write the generators with no superscript and $\G = \bO$. Here, recall that $c_{11}, c_{12}, c_2, c_3 \in \Q[X]$.
\end{rmk}

\begin{rmk}
    Our generators are related to the generators $v_1$, $v_2$, and $v_3$ introduced in~\cite{yy04} by the formulae
    \[ v_1 = -\msc{E}_1, \qquad v_2 = -\msc{E}_2, \qquad \text{and} \qquad v_3 = \msc{E}_3 - \msc{E}_2 X. \]
\end{rmk}

\begin{lem}
    The ring $\tilde{\msc{R}}$ is closed under the derivative $D$.
\end{lem}

\begin{proof}
    A direct computation yields
    \begin{align*}
        D \msc{E}_1 &= -X(\msc{E}_1 + r_0) - \msc{E}_1^2 + 2 \msc{E}_2, \\
        D \msc{E}_2 &= -X \msc{E}_2 - \msc{E}_1 \msc{E}_2 + \msc{E}_3, \\
        D \msc{E}_3 &= r_1 X - X \msc{E}_3 - \msc{E}_2^2.
    \end{align*}
\end{proof}

\begin{thm}[Reduction of generators]\label{thm:reduction}
    Let $g > 1$. Then $P_g \in \tilde{\msc{R}}$.
\end{thm}

\begin{proof}
    First, note that by definition, we have $P_g = \tilde{P}_g$. We will now prove that all $\tilde{P}_{g,m} \in \tilde{\msc{R}}$ by induction on the lexicographic order in $(g,m)$. Recall that $f_{g,m}^{\bB}$ can be computed from $\tilde{P}^{\bB}_{h \leq g,m,n}$ by the geometric quantization of $\msc{R}^{\bB}$. The contribution of each stable graph to this quantization is given by the following construction:
    \begin{itemize}
        \item At each leg, we place $\varphi_1$ or $\varphi_0 \psi$;
        \item At every edge, we place the bivector
            \begin{align*}
                &\tilde{E}_{\varphi\varphi} (\varphi_1 \otimes \varphi_1) + \tilde{E}_{\varphi\psi}(\varphi_1 \otimes \varphi_0 \psi + \varphi_0 \psi \otimes \varphi_1) + \tilde{E}_{\psi\psi}(\varphi_0\psi \otimes \varphi_0\psi) \\
                =&\msc{E}_1 (\varphi_1 \otimes \varphi_1) + \msc{E}_2(\varphi_1 \otimes \varphi_0 \psi + \varphi_0 \psi \otimes \varphi_1) + \msc{E}_3(\varphi_0\psi \otimes \varphi_0\psi);
            \end{align*}
         \item At every vertex, we place the linear map $\varphi_1^{\otimes m} \otimes (\varphi_0\psi)^{\otimes n} \mapsto \tilde{P}_{g,m,n}$.
    \end{itemize}
    
    The base cases are $\tilde{P}_{1,0,1} = \frac{\chi(Z)}{24}-1$ and $\tilde{P}_{0,3} = 1$. The dilaton equation implies that if $\tilde{P}_{g,m} \in \tilde{\msc{R}}$, then $\tilde{P}_{g,m,n} \in \tilde{\msc{R}}$ for all $n$. Now we assume that $\tilde{P}_{h,\ell,n} \in \tilde{\msc{R}}$ for all $(h,\ell) < (g,m)$. Then we know $f_{g,m}^{\bB} \in \Q[X]$ by~\Cref{cor:bmodelfeynman}. Computing it by the stable graph sum, we see
    \[ f_{g,m}^{\bB} = \tilde{P}_{g,m} + \sum_{\Gamma \text{ non-leading}} \frac{1}{\ab|\Aut \Gamma|} \on{Cont}_{\Gamma}. \]
    By the inductive hypothesis and the formula for the edge contributions, $\on{Cont}_{\Gamma} \in \tilde{\msc{R}}$ for any non-leading $\Gamma$. The desired result follows immediately.
\end{proof}

\begin{proof}[Proof of~\Cref{thm:hae}]
    The second equation~\Cref{eqn:hae2} is equivalent to~\Cref{thm:reduction} by the results of~\cite{yy04}, so we only need to prove~\Cref{eqn:hae1}.
    We proceed by differentiating the quantization action. By definition, we have
    \[ e^{P^{\bB}(\hbar, x, y-E_{\psi}x)} = e^{-\hbar V^{\bB}(\partial_{\bt}, \partial_{\bt})} e^{f^{\bB}(\hbar, x, y)}. \]
    Applying $\partial$ being either $\partial_A$, $\partial_B$, $\partial_{B_2}$, or $\partial_{B_3}$, we see that
    \begin{align*}
        e^{\tilde{P}^{\bB}(\hbar, x, y)} &= -\hbar \partial V^{\bB}(\partial_{\bt}, \partial_{\bt}) e^{V^{\bB}(\partial_{\bt}, \partial_{\bt})} e^{f^{\bB}(\hbar, x, y)} \\
        &= -\hbar \partial V^{\bB}(\partial_{\bt}, \partial_{\bt}) e^{\tilde{P}^{\bB}(\hbar, x, y)}.
    \end{align*}
    Making the change of variables $\bt' = (x',y') \coloneqq (x, y-E_{\psi}x)$, we then see that
    \[ V^{\bB}(\partial_{\bt}, \partial_{\bt}) = \frac{1}{2} E_{\varphi\varphi} \pdv[order=2]{}{x'} + E_{\varphi\psi} \pdv{}{x',y'} + E_{\psi\psi} \pdv[order=2]{}{y'} \eqqcolon V^{\bB,\ms{small}}(\partial_{\bt'}, \partial_{\bt'}). \]
    From now on, we will replace $x'$ by $x$ and $y'$ by $y$ for simplicity. We now see that
    \[ e^{P^{\bB}(\hbar, x, y)} = \partial V^{\bB, \ms{small}}(\partial_{\bt}, \partial_{\bt}) e^{P^{\bB}(\hbar, x, y)}. \]

    First, $\partial_A V^{\bB, \ms{small}}(\partial_{\bt}, \partial_{\bt}) = \frac{1}{2} \pdv[order=2]{}{x}$, so we obtain
    \[ \partial_A P^{\bB}(\hbar, x, y) = -\frac{\hbar}{2} \pdv[order=2]{}{x} P^{\bB}(\hbar, x, y) - \frac{\hbar}{2} \ab(\pdv{}{x} P^{\bB}(\hbar, x, y))^2. \]
    Setting $x=y=0$, we see that
    \[ \sum_{g \geq 2} \hbar^{g-1} \partial_A P_g = -\frac{1}{2} \ab(\sum_{g \geq 2} \hbar^{g-1} P_{g-1,2} + \ab(\sum_{g_1,g_2>0} \hbar^{g_1+g_2-1} P_{g_1} P_{g_2})). \]
    Taking the coefficient of $\hbar^{g-1}$ on both sides, we obtain~\Cref{eqn:hae1}.
\end{proof}

\printbibliography

\end{document}